\documentclass{amsart}
\usepackage{amsxtra,amssymb,amsthm,amsmath,latexsym, epsf}

\theoremstyle{plain}

\newtheorem*{theorem}{Theorem}

\def\nd{\noindent}
\def\ds{\displaystyle}

\def\C{{\mathbb C}}

\def\Z{{\mathbb Z}}

\def\Im{\hbox{\ Im\ }}
\def\Re{\hbox{\ Re\ }}
\def\oN{\buildrel\circ\over N}
\def\~x{\widetilde\times}
\def\bysame{\rule{.5in}{.005in}}
\def\cop{\bot\hskip-.075in\bot}

\begin{document}

\title[~]
{
Families of four dimensional manifolds that
become mutually\\ diffeomorphic after one stabilization
}

\author{  David Auckly }
\address{
Department of Mathematics \\
Kansas State University\\
        Manhattan, KS  66506-2602, USA \\
        dav@math.ksu.edu}

\thanks {This work was partially supported by grant CMC 9813183
from the National Science Foundation.  }


\date{}

\maketitle\thispagestyle{empty}





It is well known that two homotopy equivalent, simply connected 4-manifolds
become diffeomorphic after taking the connected sum with enough copies of
$S^2\~x S^2$ [20]. The same result is true with $S^2\~x S^2$ replaced by 
$S^2\times S^2$, and similar results are known for special families of 4-manifolds when $S^2\~x S^2$ is replaced by other manifolds. Taking the connected sum with one of these specific manifolds is called stabilization.
For this paper, we will only consider connected sums with $S^2\~x S^2$, and
stabilization will refer to taking the connected sum with this specific
manifold. Most of the arguments in this paper can be easily modified to
address other summands as well. 
Many families of distinct homotopy equivalent simply connected 4-manifolds
that become mutually diffeomorphic after one stabilization are known [15]. 
There is, in fact, no known pair of homotopy equivalent simply connected 4-manifolds which are not diffeomorphic after one stabilization.

In this paper, we will introduce 
a cut and paste move, called a geometrically null log transform, and prove that any two manifolds related by a sequence of these moves become diffeomorphic after one stabilization. To motivate the cut and paste move, we will use the symplectic fiber sum,
and a construction of Fintushel and Stern
to construct several large families of 4-manifolds.
We will then proceed to prove that the members of any one of these families
become diffeomorphic after one stabilization. Finally, we will
compute the Seiberg-Witten invariants of each member of each of the
families.

Even though the  Donaldson and Seiberg-Witten invariants
can distinguish some homotopy equivalent four-manifolds,
these invariants cannot directly distinguish manifolds of the form
$X\# S^2\~x S^2$. This is because both invariants
are trivial on 4-manifolds with an $S^2\~x S^2$ summand,
provided that the second positive betti number of the remaining summand is 
positive [19].
Apriori, it is possible that
$X\#S^2\~x S^2\cong Y\#S^2
\~x S^2$ implies some relation between
the Seiberg-Witten invariants of $X$ and the Seiberg-Witten
invariants of $Y$. The first reason for considering a specific set of families in this paper is to show that no
simple relation between Seiberg-Witten invariants is implied by equivalence
after one stabilization.

If it was known that any pair of homotopy equivalent simply connected 4-manifolds are related by a sequence of geometrically null log transforms,
it would follow that any two such manifolds become equivalent after one
stabilization. It is known that any manifold homotopy equivalent to a 
simply connected 4-manifold may be constructed by removing a contractible
4-manifold and reglueing it via an involution [3], [14]. This motivates the question: Is it possible to modify the proof of the decomposition theorem
to find a finite set of moves which could be used to pass between any two
homotopy equivalent 4-manifolds? The contractible piece is known as a cork. 
A second reason for constructing specific families is to study the effect that applying a geometrically null log transform to one manifold of 
a pair of homotopy equivalent simply connected 4-manifolds has on the cork.  

I would like to thank Bob Gompf for
a helpful conversation regarding this material.

\section*{Families of 4-Manifolds.} 

All of the 4-manifolds explicitly considered in this paper are formed by applying
a cut and paste operation, the fiber sum,
to copies of a standard building block, called the K3 surface. This section
begins with a short description of the K3 surface. (See the book by Harer, 
Kas, and Kirby for more information about the K3 surface [9].)
This section will end with explicit handle decompositions of the
4-manifolds contained in the specific families considered in this paper.

Recall that the K3 surface is essentially the quotient of a
4-torus by an involution.
The group, $\Z_2$ acts on $T^4$ via the map:
$$ \varepsilon:T^4=\frac{\C^2}{\Z[i]^2}\to
   T^4;\quad \varepsilon ([x,y])=[-x,-y].
   $$
It also acts on $\C P^2$ via
$$ \eta: \C P^2\to \C P^2;
   \quad \eta([x:y:z])=[-x:-y:z].   $$
There are 16 fixed points on $T^4$, namely,
 $ \left(1/2\Z[i]\right)^2 \big/ \left(\Z[i]\right)^2$,
and the fixed point set in $\C P^2$ is
 $ \{[0:0:1]\} \cup \{[x:y:0]\}$.
We may cut invariant neighborhoods of the 16 fixed points out of
$T^4$ and glue in 16 copies of the complement of an invariant
neighborhood of $\{[0:0:1]\}\subseteq \C P^2$,
to get a $\Z_2$ action on
$T^4\#\overline{\C P^2}^{\#16}$.
The bar refers to the fact that $\C P^2$ is taken with the
opposite orientation.
The quotient of
$T^4\#\overline{\C P^2}^{\#16}$ by $\Z_2$ is 
the $K3$ surface. It is manifold essentially because the
quotient of the disk,
$$ D_{[x_0:y_0]} = \{ [x:y:z] \vert [x:y] = [x_0:y_0]
   \ \&\  |z|^2\leq |x|^2 +|y|^2 \},
   $$
in $\C P^2$ is also a disk.

All of the examples that we construct will be obtained by cut and paste
along three tori in the $K3$ surface.
Let
$$ \begin{array}{rl}
   T_1 & =\left\{ [(x,1/3+1/3i)] \in K3 \vert x\in \C \right\},\cr
   T_2 & = \left\{ [(x,y)] \in K3 \vert \Im x=\Im y=1/4 \right\}\cr
   \text{and}\quad T_3
      & = \left\{ [(x,y)] \in K3\vert \Im x=\Re y=1/5 \right\}.
      \end{array} $$
Let $X_N$ be the manifold obtained by fiber summing $N$ copies of the
$K3$ surface  together along $T_3$ in one copy and $T_1$ in the
next copy.
 (See figure 1)

\vskip -.1in
\begin{figure}[htbp]
  \begin{center}
    \leavevmode
    \epsfysize=.8in
    \epsfbox{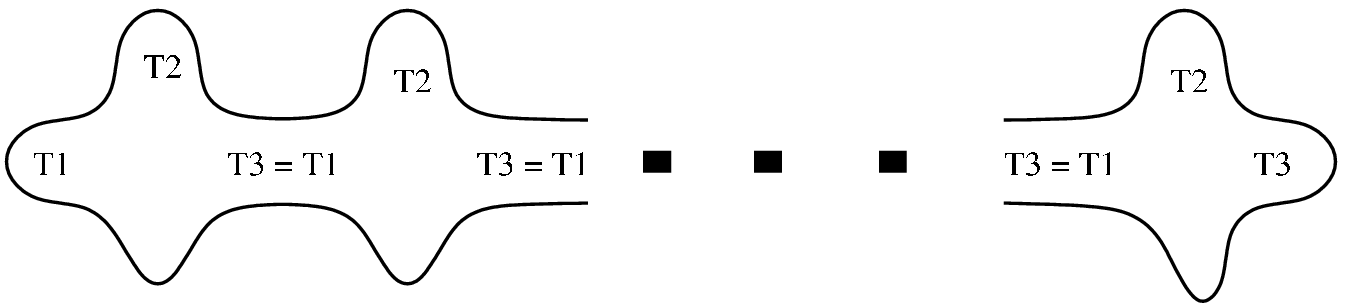}
  \end{center}
\end{figure}
\vskip-.1in

\centerline {\bf Figure 1: $X_N$.} 

\vskip.15in

\nd The fiber sum of $(X,S)$ and $(Y,T)$  is
$\big(X-\buildrel\circ\over{N} (S)\big)\cup_{\partial N(T) =\partial N(S)}
\big(Y-\buildrel\circ\over{N} (T)\big)$.
It will be denoted by $(X,S)\#(Y,T)$.
If $S$ and $T$ are symplectic submanifolds with opposite
self-intersection numbers, the fiber sum will also be symplectic [8]. 
The definition of the fiber sum requires an orientation reversing glueing map from the boundary
of a tubluar neighborhood of $S$ to the boundary of a tubular 
neighborhood of $T$. Every thing that we will assert about the manifolds,
$X_N$ will be independent of the glueing maps. To be definite one could
choose $\varphi : \partial N(S) \to \partial N(T)$ given by, 
$\varphi (x_1+1/5 i + 10^{-2}i\cos(\theta), 1/5 +10^{-2}\sin(\theta)+x_2i)\,=\,
(x_1+x_2i, 1/3+10^{-2}\cos(\theta)+1/3i-10^{-2}i\sin(\theta))$.
The manifold $X_N$  has $N+2$ of tori
the $T_i$ remaining.
Copies of $S^1\times S^3$ may be fiber summed onto these
remaining tori, each along $S^1$ cross a knot (take the glueing map which
identifies the 0-framed longitude of the knot with a meridian or the torus). Fintushel and Stern
proved a remarkable formula relating the Alexander polynomial of a knot
to the change in the Seiberg-Witten invariant of a manifold after
fiber summing with $S^1\times S^3$ along $S^1$ cross the knot. [5].
This formula will be used to compute the Seiberg-Witten invariants at
the conclusion of this paper.
All of the manifolds obtained from a fixed $X_N$, by fiber summing
with $S^3\times S^1$ as above are homotopy equivalent. We will show that
all members of the family of manifolds obtained from a fixed $X_N$ become
diffeomorphic after one stabilization. The last section of the paper
describes the Seiberg-Witten invariants of these manifolds.

A well-known handle decomposition of the K3 surface is given
in the book by Harer, Kas, and Kirby [9].
This handle decomposition has 24 handles, the minimal number
of handles in a handle decomposition of the K3 surface.
Other 4-manifolds will require even more handles. Because of this complexity,
it is useful to decompose 4-manifolds into a union of compact pieces and
then describe handle decompositions of the pieces. One important  piece
of the K3 surface is the Gompf Nucleus. By definition, this is
a neighborhood of the union of a cusp fiber and a section [6].
The nucleus of K3 will be denoted by $N_2$.
It may be constructed by attaching three two-handles to $T^2\times D^2$
  (see figure 2).


\vskip -.25in
  \begin{center}
    \leavevmode
    \epsfysize=1.9in
    \epsfbox{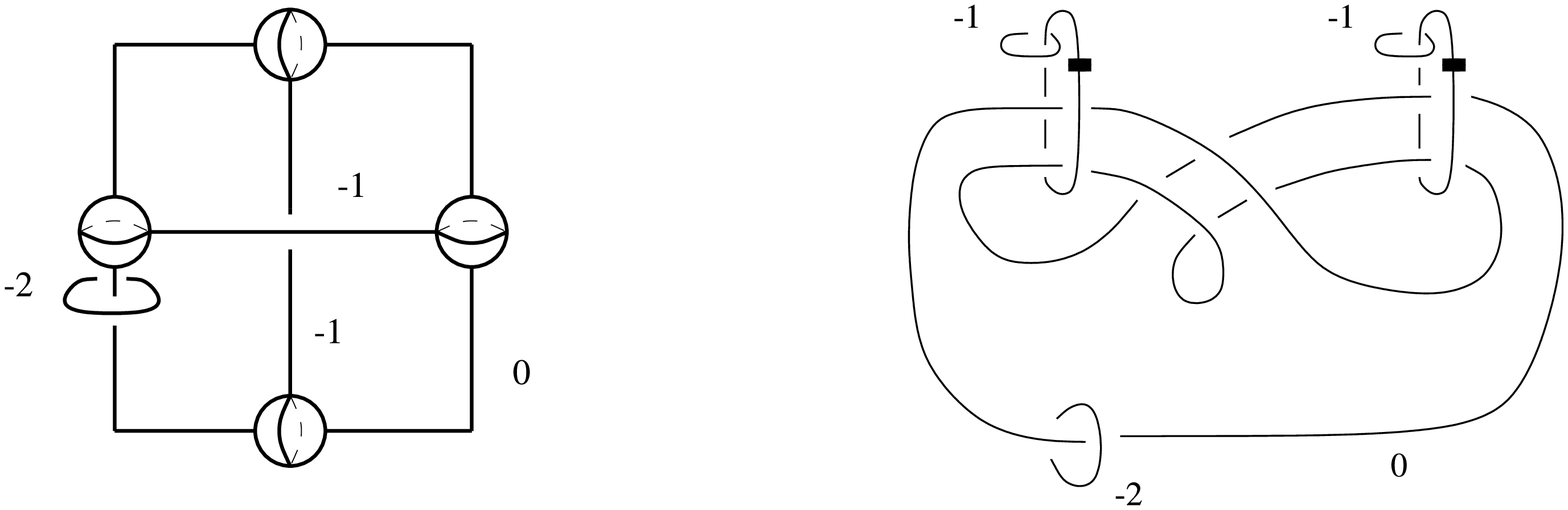}
  \end{center}
\vskip-.1in
   \centerline{\bf Figure 2:  $N_2$.}
\vskip.2in



There are three disjoint copies of the nucleus in the K3 surface.
Each one contains one of the, $T_i$, tori described above as $T^2\times\{0\}$
in figure 2. Given a handle decomposition of a 4-manifold
with boundary, it will be useful to denote a collar of the boundary
by putting an $I$ on each handle.
For example, figure three displays handle decompositions of
$N(\partial(T^2\times D^2))$ and $N_2-\oN(T^2)$. 

\begin{figure}[htbp]
  \begin{center}
    \leavevmode
    \epsfxsize=5.5in
    \epsfbox{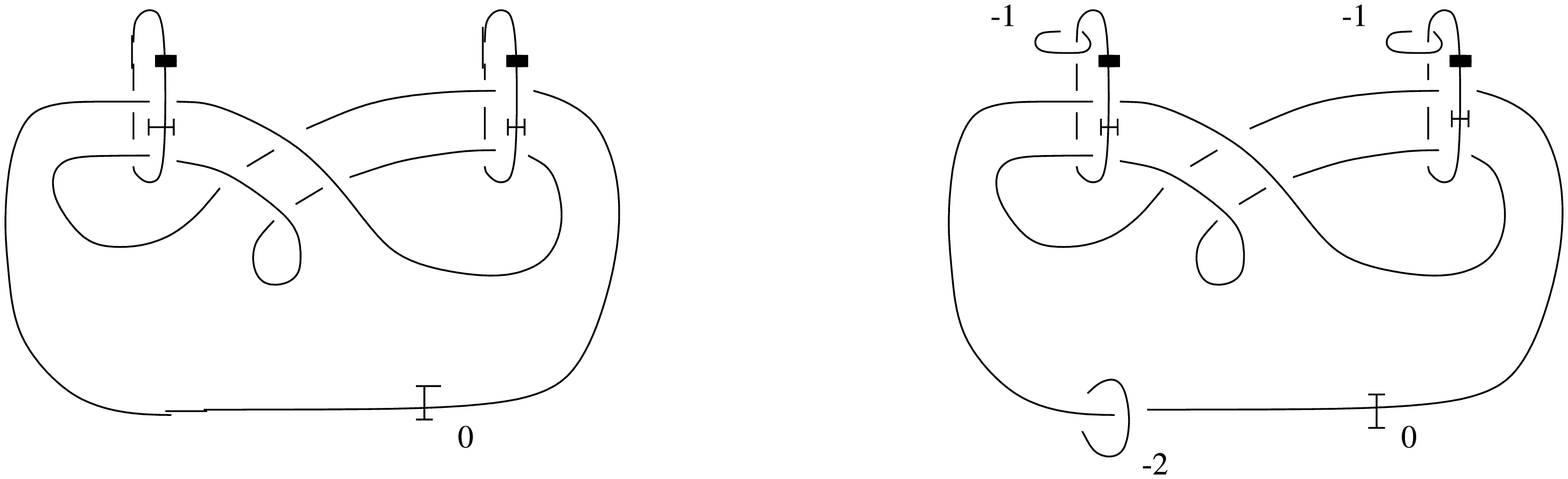}
  \end{center}
   \centerline{\bf Figure 3:  $N(\partial (T^2\times D^2))$,
  $N_2-\oN(T^2)$.}
\end{figure}

To construct a handle decomposition of the fiber sum of a pair of
nuclei, we will turn a copy of
$N_2-\oN(T^2)$ upside down and glue it to a second copy of 
$N_2-\oN(T^2)$.
To turn a handle decomposition upside down, first reverse the orientation
(reverse every crossing and framing), then double.
Assuming that the original manifold has no 3-handles,
attach one 0-framed 2-handle to the co-core of
each original 2-handle, then delete the original manifold
(add $I$'s to all of the original components). Figure 4 displays
$N_2-\oN(T^2)$ turned upside down and a fiber sum of a pair of nuclei
constructed by glueing the \linebreak $N_2-\oN(T^2)$ from figure 3 to the $N_2-\oN(T^2)$
from figure 4.

Turn now to the construction of handle decompositions for manifolds
of the form $M^3\times S^1$.
Restrict the boundary of $M$ to be a disjoint union of tori. If $M$ is
described by surgery to the complement of a link is $S^3$,
there will be two approaches for constructing handle decompositions
for $M^3\times S^1$.
Both methods begin by constructing a handle decomposition for
$M^3\times I$. The first method is to
pick a tunnel system for the link
$L=L_1\cop L_2$ when  $M$ is obtained by Dehn filling on $L_1$.
This tunnel system may be used to construct a handle decomposition
of $S^3-\oN(L_1\cop L_2)$.
This is easily translated into a handle decomposition of
$M^3$ and then $M^3\times I$ (see p.250 of [18] for this process applied to
the Poincar\'e homology sphere.)


  \begin{center}
    \leavevmode
    \epsfxsize=6in
    \epsfbox{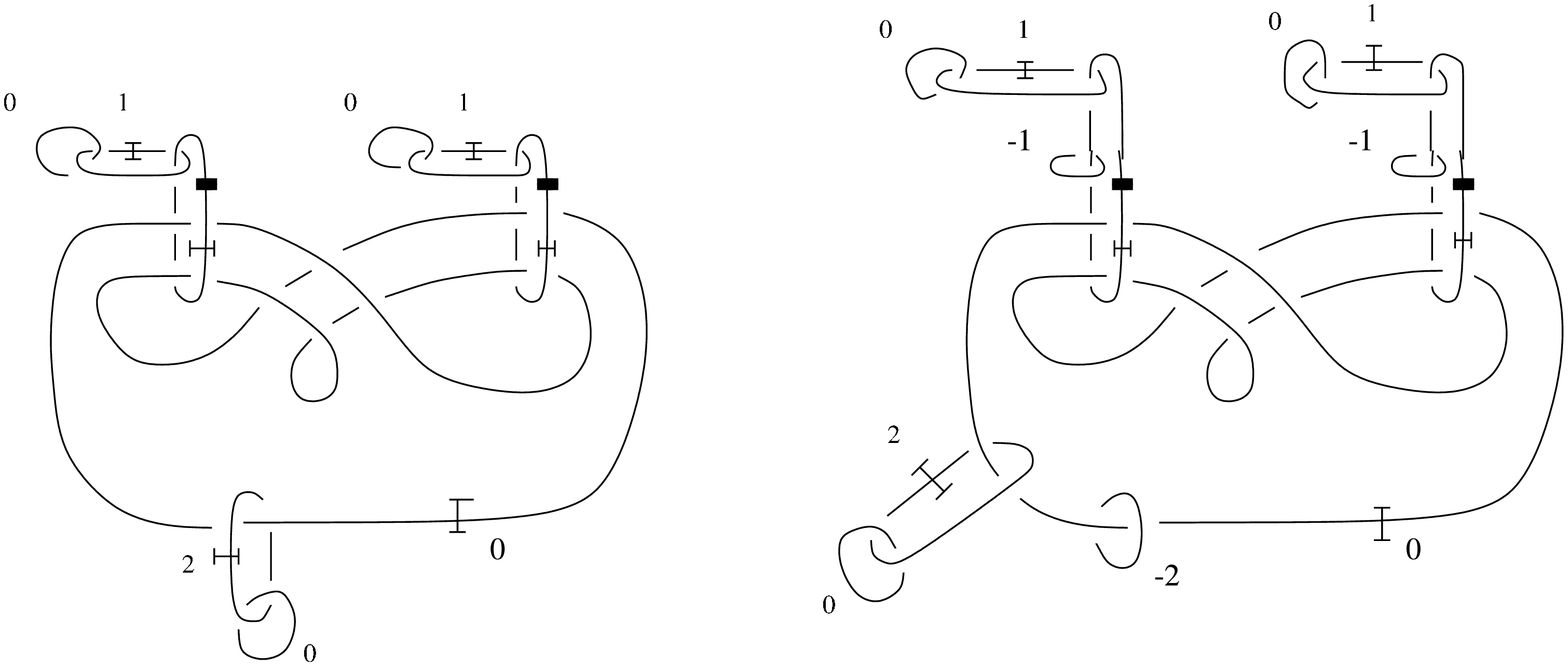}
  \end{center}
   \centerline{\bf Figure 4: $N_2-\oN(T^2),(N_2,T^2)\#(N_2,T^2)$,}

\vskip .1in


The second approach is based on the observation that proves that $K\#-K$ is
slice for any knot, $K$ [1], [2]. Namely, $(K-\oN(pt))\times I$ is a slice disk for $K\#-K$. For any link, $L$, $I\times (S^3-L)$ may be described as the
exterior of a surface, $F$, in $D^4$. The surface, $F$, is constructed in the same way as the slice disk for $K\#-K$. If $M^3$ is surgery on $L$, a 
handle decomposition of $I\times (S^3-L)$ may easily be converted into 
a decomposition of $M^3 \times I$.
To begin the description of $D^4-\oN(F)$, notice that
$$I\times S^3-\oN(I\times L)=I\times S^3-\oN(I\times pt))-
\oN(I\times (L-\oN(pt))=D^4-\oN(F).$$

\nd Figure 5 shows a typical link and the frames of the movie obtained by
intersecting $D^3\times\{t\}$ with
the canonical cobordism in $D^3\times I=D^4$.


\begin{figure}[htbp]
  \begin{center}
    \leavevmode
    \epsfysize=3in
    \epsfbox{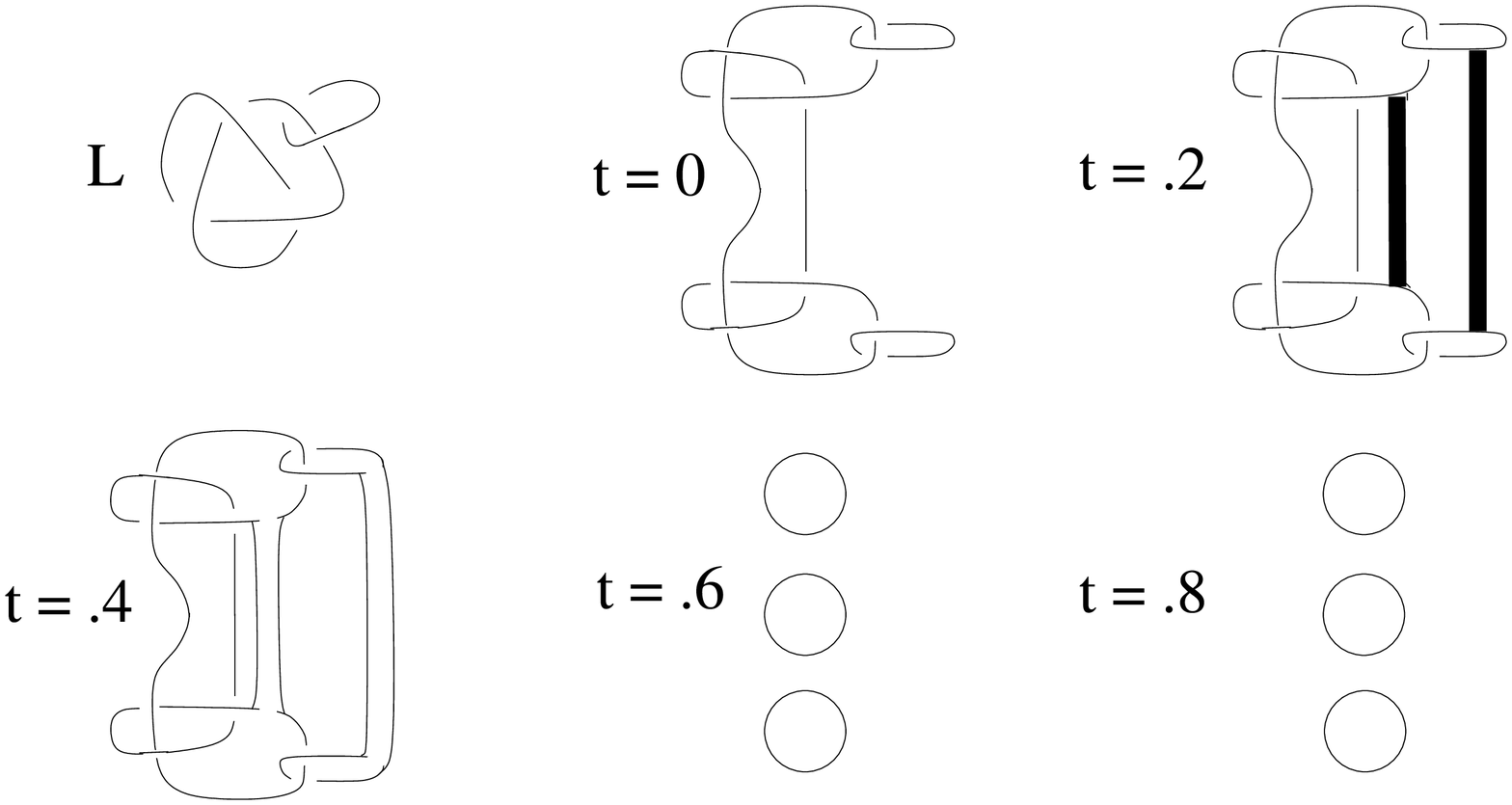}
  \end{center}
   \centerline{\bf Figure 5: $L$ and $F$.}
\end{figure}


\nd Note that
\begin{equation}
  \begin{aligned}
  (D^4-\oN(F))\cap(D^3\times[.6,1])
   &=D^4-\oN((D^2)^{\cop k})=(D^2-\oN(k\, pts))\times D^2
  \notag\\
   &=(D^2 \ds\operatorname*{\cup}_{k(\partial \,D^2)\times D^1}
       k\,D^1\times D^1)\times D^2
  \notag\\
   &=D^4\ds\operatorname*{\cup}_{k(\partial\,D^1)\times D^3} kD^1\times D^3.
  \notag\end{aligned}\end{equation}
This will allow us to describe
a handle decomposition for $D^4-\oN(F)$
in terms of a handle decomposition for $F$.
So far we see that 0-handles in $F$ correspond to 1-handles in
$D^4-\oN(F)$ (this is the previous computation). Figure 6 displays a neighborhood of a 1-handle
in $F$ embedded in $D^4$. The cylinder around the band is a 2-handle
in $D^4-\oN(F)$.

  \begin{center}
    \leavevmode
    \epsfxsize=6in
    \epsfbox{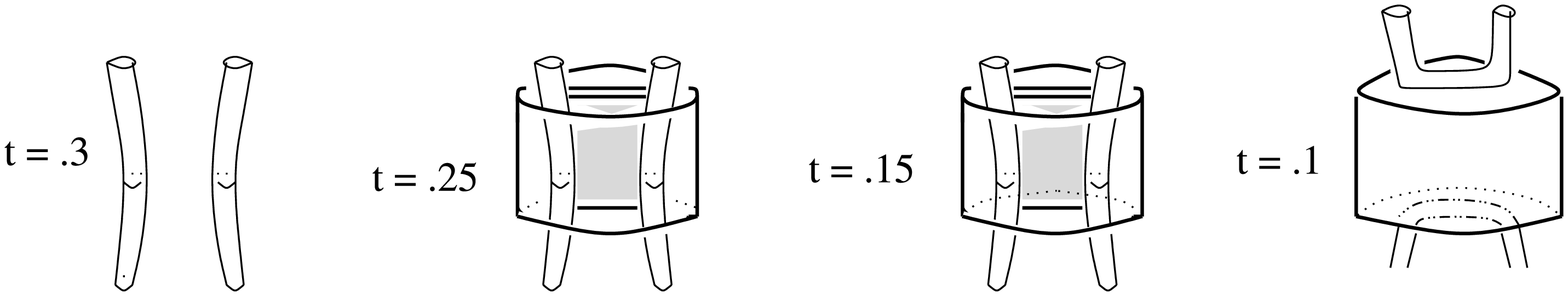}
  \end{center}
   \centerline{\bf Figure 6: Neighborhood of a 1-handle.}


\nd This illustrates the fact that 1-handles in $F$ correspond to
2-handles in $D^4-\oN(F)$. It also enables one to construct handle
decompositions for $(S^3-\oN(L))\times I$.
Dehn filling is accomplished by attaching a 2-handle and then
attaching a 3-handle. This will complete a handle decomposition
of $M^3\times I$.
The special cases when $M$ is $D^3$, or $S^1\times D^2$, or $S^2\times D^1$
are instructive, when extending a handle decomposition of $M^3\times I$
to a handle decomposition of $M^3\times S^1$.
In general, a $(k+1)$-handle is added for every $k$-handle of $M^3\times I$.

We can apply these ideas to $M=S^3-\oN(K)$.
Let the knot $K$ be expressed as the closure of a braid, $\beta$,
in such a way that the black board framing of $K$ is the zero framing.
The result is the handle decomposition for
$(S^3-\oN(K))\times S^1$ displayed in figure 7.
Figure 7 also has a handle decomposition of
$(N_2, T)\#(S^3\times S^1, K\times S^1)$ obtained from
$(S^3-\oN(K))\times S^1$ by gluing on an
$N_2-\oN(T)$.

\begin{figure}[htbp]
  \begin{center}
    \leavevmode
    \epsfxsize=6in
    \epsfbox{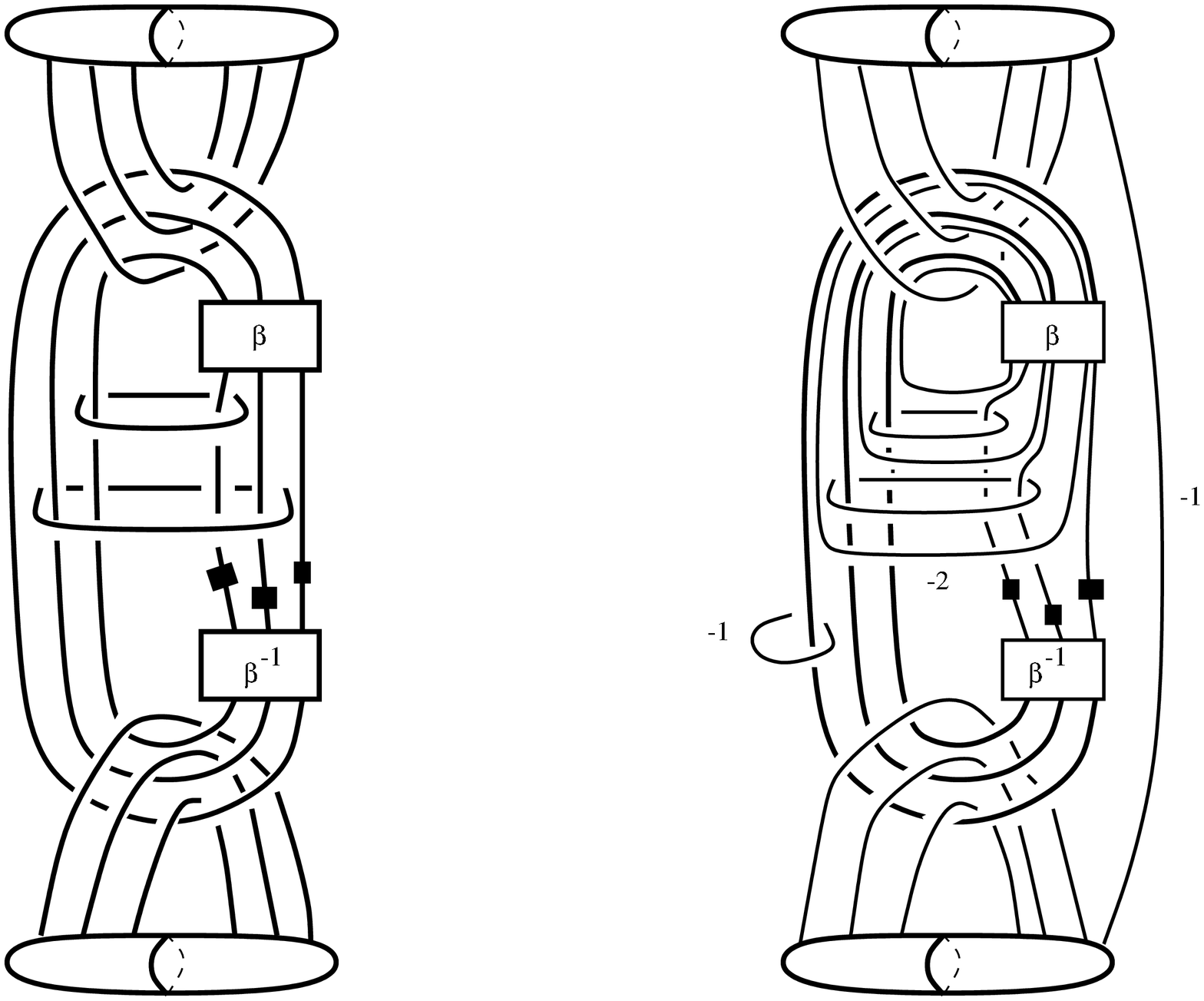}
  \end{center}
   \centerline{\bf Figure 7:
$(S^3-\oN(K))\times S^1$ and $(N_2,T)\#(S^3\times S^1, K\times S^1)$.}
\end{figure}


\begin{figure}[htbp]
  \begin{center}
    \leavevmode
    \epsfxsize=6in
    \epsfbox{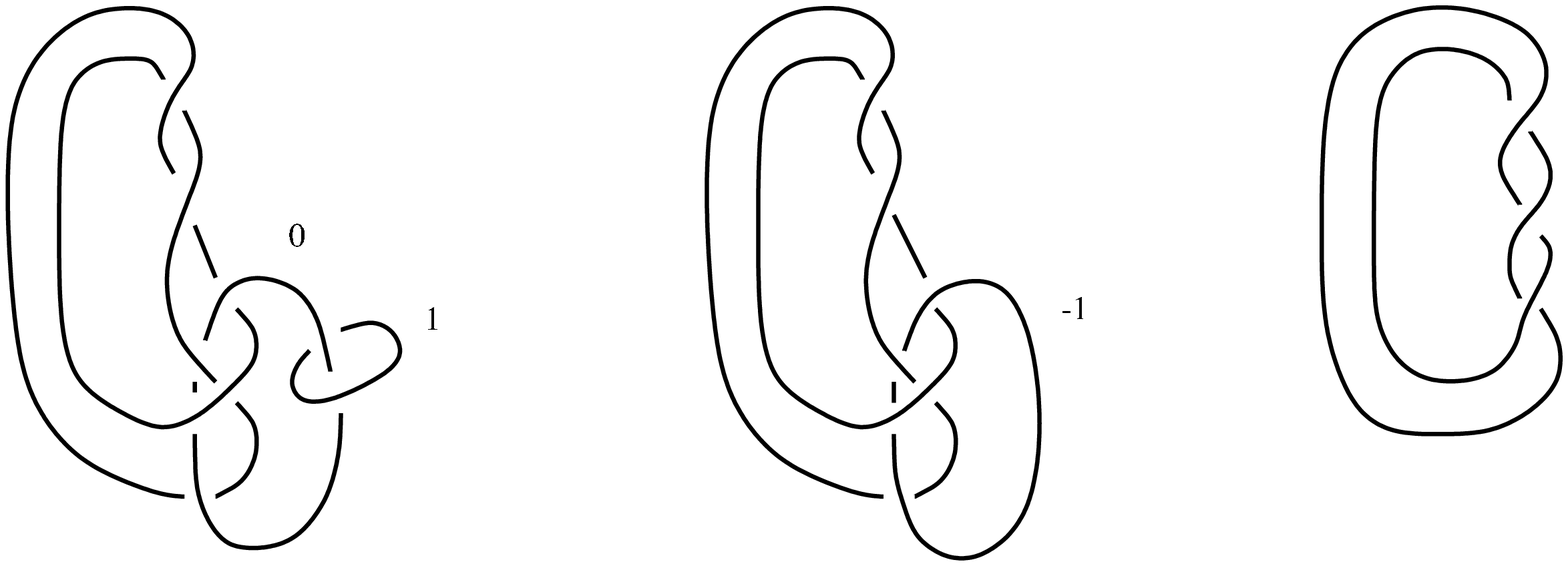}
  \end{center}
   \centerline{\bf Figure 8:
Different descriptions of the same manifold.}
\end{figure}

There are many different surgery descriptions of any given 3-manifold
(see figure 8).
Any of these descriptions will produce a handle decomposition of
$M^3\times S^1$.
It is an interesting exercise to see how various 3-manifold moves translate
into sequences of  handle slides and handle pair birth/deaths.
In particular, it is interesting to see how Markov
moves on the braid, handle slides, and Kirby moves effect the
4-dimensional handle decomposition.

Notice that any knot can be converted to the unknot by a sequence
of $\pm 1$ surgeries. This will enable us to
understand the fiber sum with $S^3\times S^1$ along a
complicated knot crossed with the circle using one simple move.
We will come back to this later in this paper.


\section*{Stabilization}

For this paper, stabilizing a 4-manifold will simply
refer to  taking the connected
sum with $S^2\~x S^2$. The manifold, $S^2\~x S^2$ is
the nontrivial $S^2$ bundle over $S^2$.
It may also be described as $\C P^2\#\overline{\C P^2}$.
Stabilization is closely related to the surgery corresponding to the
addition of a five dimensional 2-handle. This surgery
amounts to replacing an $S^2\times D^3$ by a $D^2\times S^2$
in the 4-manifold. If $S^1\times\{0\}$ is homotopically trivial,
we may assume that it is contained in a 4-disk.
Since surgery on a trivial loop in the 4-disk either produces a
punctured $S^2\times S^2$ or $S^2\~x S^2$, it follows
that surgery on a null homotopic loop is the same as taking the connected
sum with either $S^2\times S^2$ or $S^2\~x S^2$ (see figure 9).

\vskip-.2in

  \begin{center}
    \leavevmode
    \epsfxsize=6in
    \epsfbox{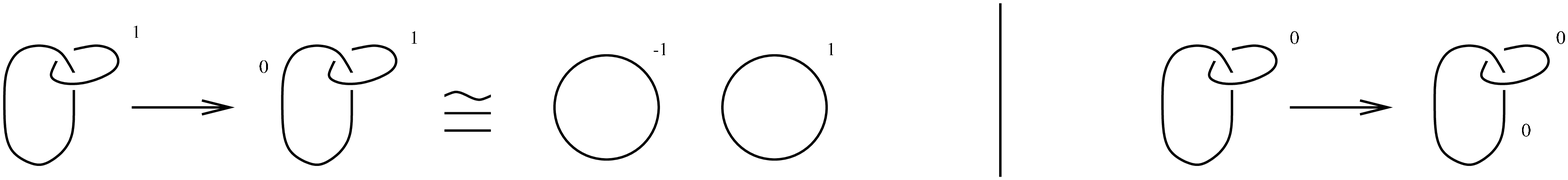}
  \end{center}
\vskip-.08in
   \centerline{\bf Figure 9: Surgery and stabilization.}

\vskip .15in

\nd Combining this with the observations that
$(S^2\~x S^2)\#(S^2\times S^2)\cong
(S^2\~x S^2)\# (S^2\~x S^2)$,
and that any five dimensional $h$-cobordism may be constructed with
just 2-handles and 3-handles proves that two homotopy
equivalent, simply connected 4-manifolds become diffeomorphic
after some number of stabilizations [10], [20].

Computing the number of stabilizations required is an interesting open
problem. For every known example, one stabilization is enough.
The main argument used to prove that one stabilization is enough is a
five dimensional handle argument due to Mandelbaum [11], [12], [13].
In fact, many manifolds are known to become diffeomorphic to
$(\C P^2)^{\# n}\#(\overline{\C P^2})^{\# m}$ after taking
the connected sum with just $\C P^2$ [15].
Many related facts may be found in [7].
If $S$ and $T$ are tori in $X$ and $Y$, the basic five
dimensional argument analyzes a natural cobordism between
$X\cop Y$ and
$(X,S)\#(Y,T)$. Let $S$ have a standard handle decomposition,
$S=h^{(0)}\cup h^{(1)}_1\cup h^{(1)}_2\cup h^{(2)}$.
The natural cobordism is then
$$ \begin{array}{rl}
  W&=(I\times(X\cop Y)\cup D^1\times h^{(0)}\times D^2\\
   &\cup D^1\times h^{(1)}_1\times D^2 \cup D^1\times h^{(1)}_2\times D^2\\
   &\cup D^1\times h^{(2)}\times D^2.
   \end{array} $$
The level of $W$ after the 1-handle, $D^1\times h^{(0)}\times D^2$,
is $X\# Y$. The level after the 2-handles,
$D^1\times h^{(1)}_1\times D^2$ and $D^1\times h^{(1)}_2\times D^2$,
is $X\#Y\#(S^2\times S^2)\times (S^2\times S^2)$.
The section of the cobordism from this level to the end is obtained
by attaching a 3-handle. By turning this section upside down, we see
that it is also obtained by attaching a five dimensional 2-handle to
$(X,S)\#(Y,T)$.
The level is therefore $(X,S)\#(Y,T)\#(S^2\times S^2)$.
Thus
$X\# Y\#(S^2\times S^2)\#(S^2\times S^2)\cong (X,S)\#(Y,T)\#(S^2\times S^2)$.
In the above argument, we assumed that $X$ and $Y$ were simply
connected, and that the framings on all of the five dimensional
2-handles are arranged so that factors of $S^2\times S^2$ appear, not
factors of $S^2\~x S^2$.

Instead of checking the framings directly, we will use the five
dimensional argument as a guide for a four dimensional handle
sliding argument that
$X_N\#(S^2\~x S^2)\cong (\C P^2)^{\#4N}\#(\overline{\C P^2})^{20 N}$.
\smallskip

  \begin{center}
    \leavevmode
    \epsfxsize=6in
    \epsfbox{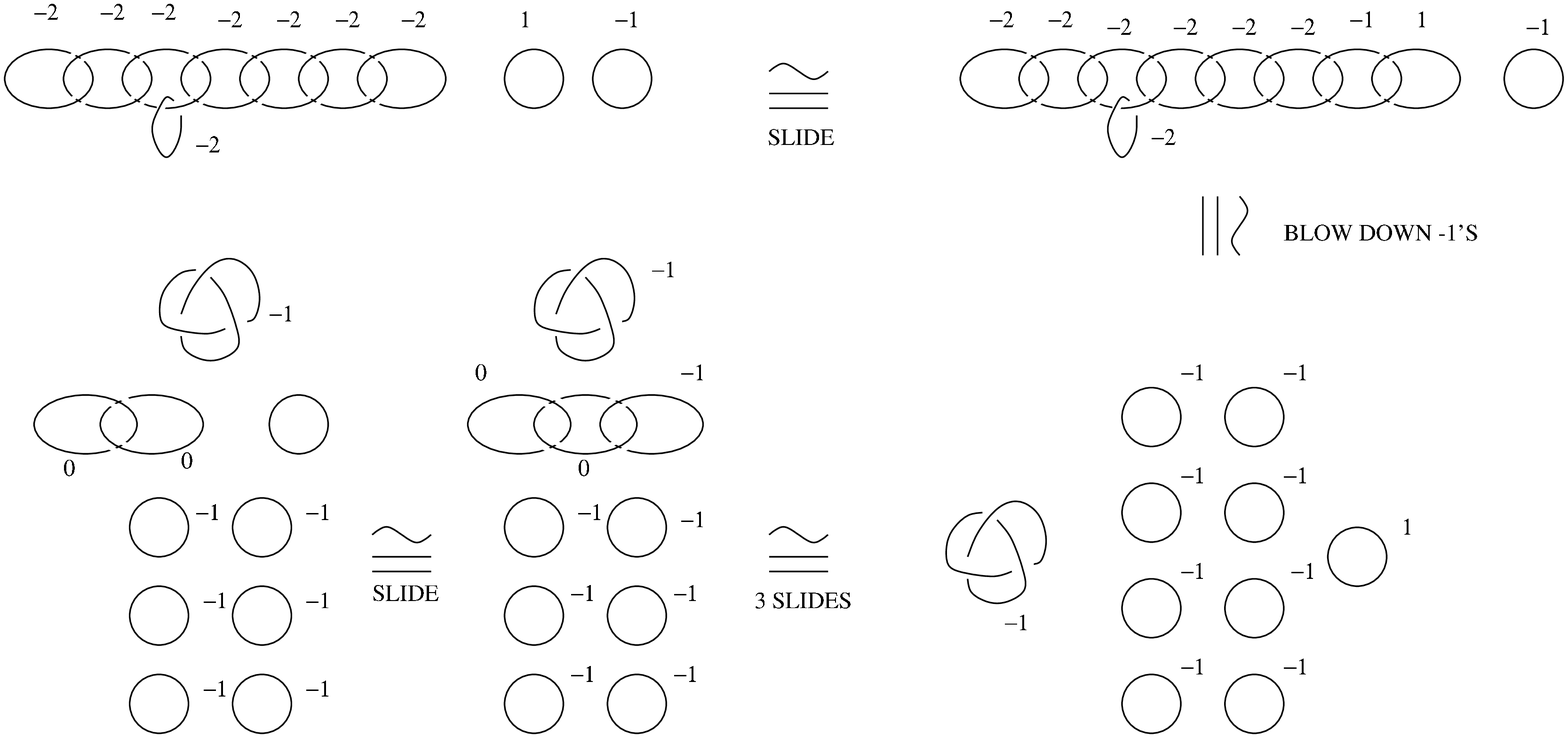}
  \end{center}
   \centerline{\bf Figure 10:  $E_8\#(S^2\~x S^2)\cong
  W_1(\overline{\C P^2})^{\#7}\#(S^2\times S^2)$.}
\vskip.1in


\nd The $E_8$ Milnor fiber is embedded in K3 disjoint from the nucleus
[9]. It follows that $E_8$ is also embedded in $X_N$ disjoint from
all of the tori used in the fiber sum. The argument begins by showing
that
$E_8\#(S^2\~x S^2)\cong W_1\#(\overline{\C P^2})^{\#7}\#(S^2\times S^2)$
(see figure 10).
Sliding the factor of $S^2\times S^2$ into the
$(N_2,T^2)\#(N_2,T^2)$
from figure 4, and performing the moves indicated in figure 11
produces figure 12.

  \begin{center}
    \leavevmode
    \epsfxsize=5in
    \epsfbox{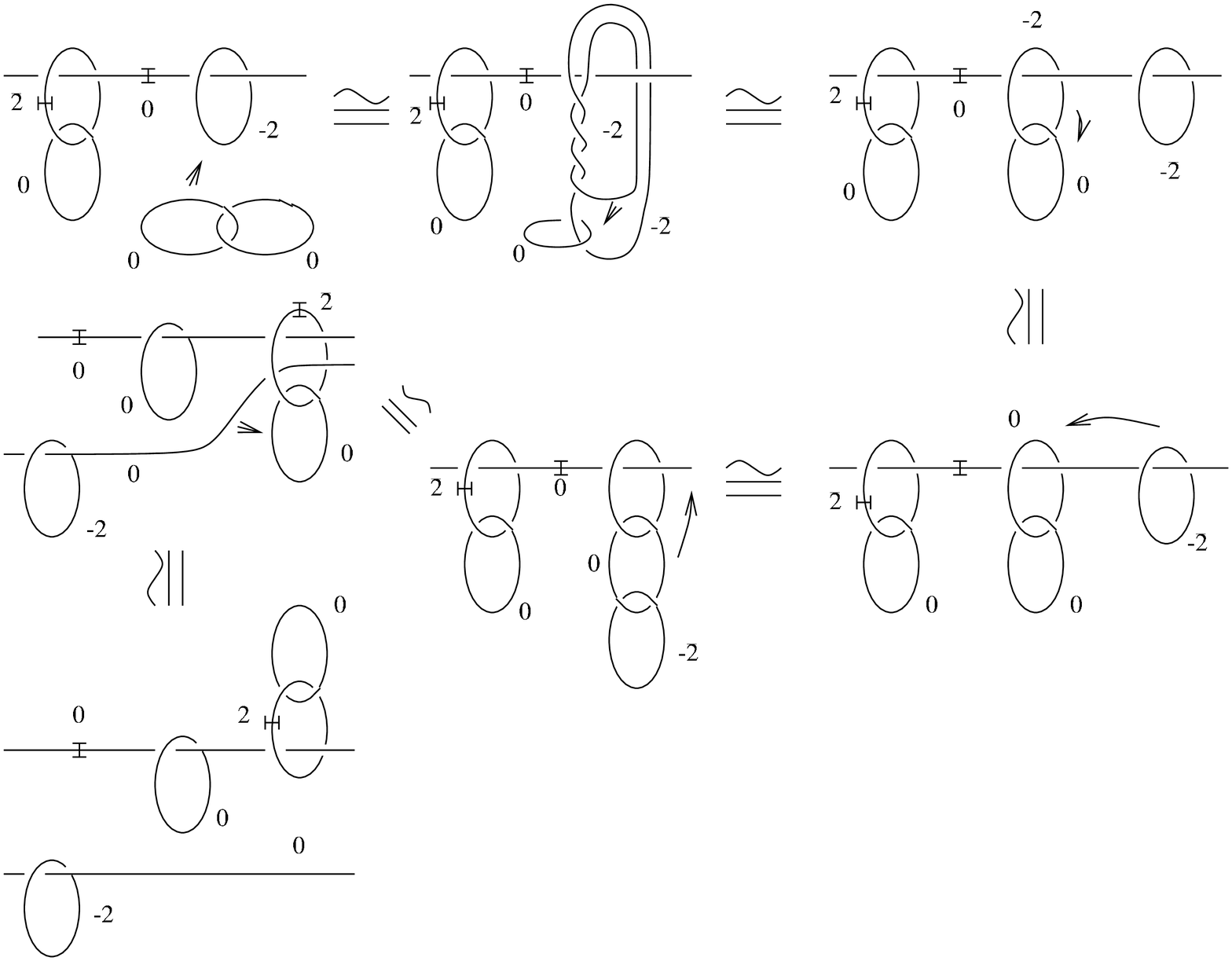}
  \end{center}
   \centerline{\bf Figure 11: Handle slides for the five
dimensional 3-handle.}
\vskip.1in



\nd The handle slides in figure 11 correspond to the last section
of the cobordism in the five dimensional argument.

  \begin{center}
    \leavevmode
    \epsfxsize=2.5in
    \epsfbox{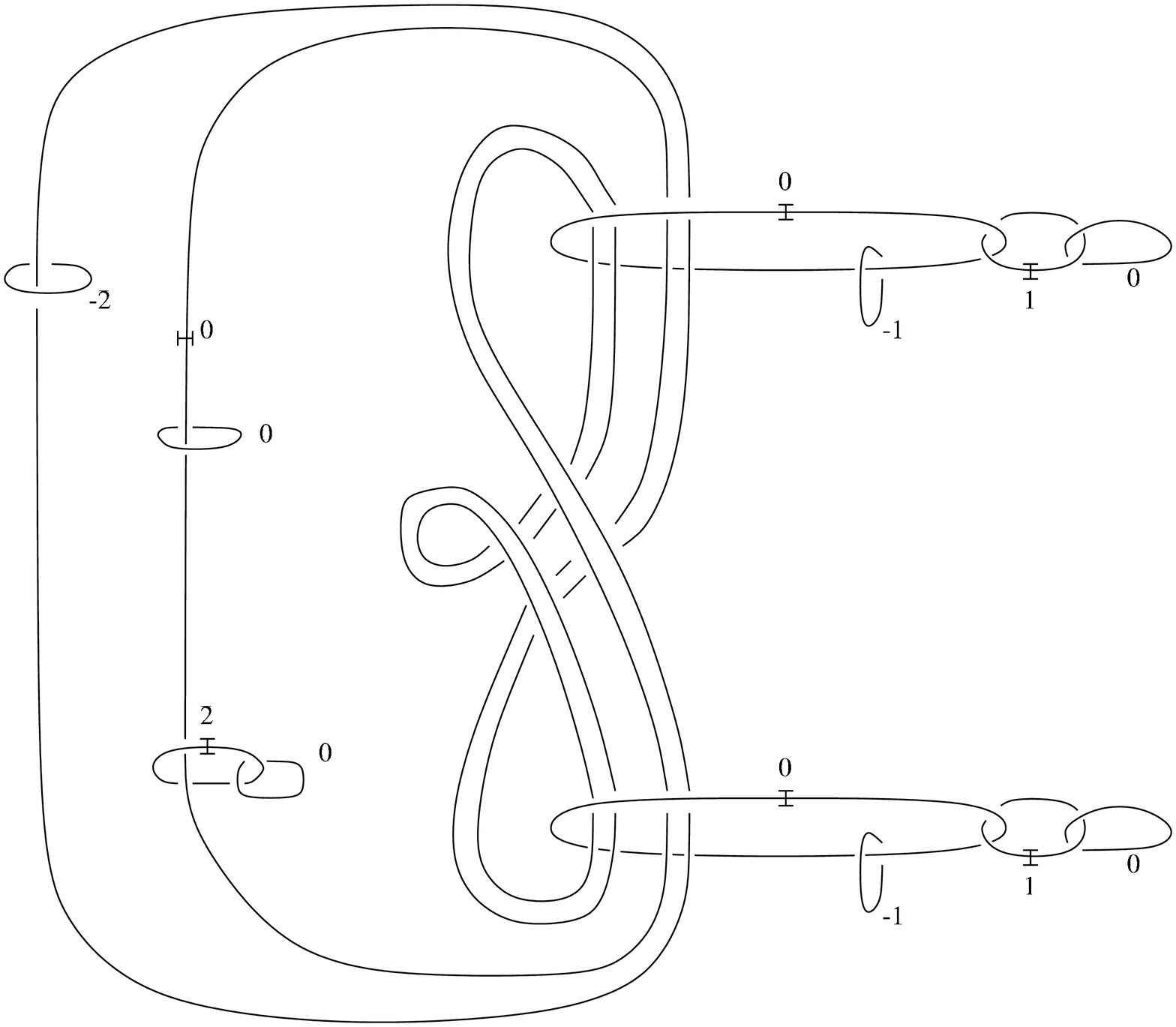}
  \end{center}
   \centerline{\bf Figure 12: $(N_2,T^2)\#(N_2,T^2)\#(S^2\times S^2)$..}

\smallskip

\nd Sliding the complicated zero framed 2-handle over the 2-handle
dual to the complicated $0I$ handle will allow the complicated zero
framed 2-handle to be pushed to the right of the figure as in figure 13.

  \begin{center}
    \leavevmode
    \epsfxsize=4in
    \epsfbox{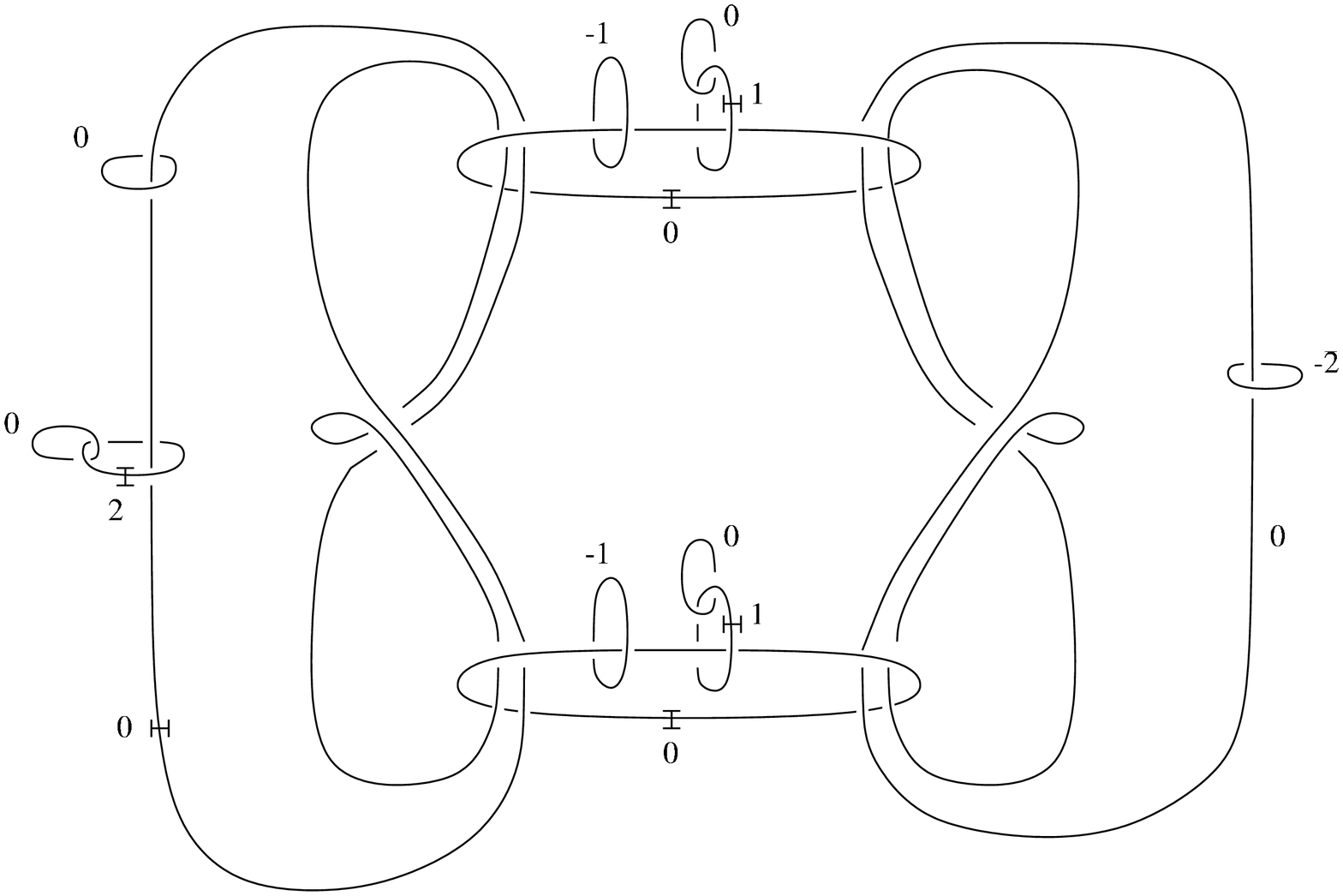}
  \end{center}
   \centerline{\bf Figure 13: $(N_2,T^2)\#(N_2,T^2)\#(S^2\times S^2)$.}

\medskip \bigskip


\nd The next step is to add two canceling 1-handle/2-handle pairs
to produce the 1-handles in the right-side-up $N_2$. This is
done in figure 14, resulting in the handle decomposition in figure 15.
\bigskip \medskip

  \begin{center}
    \leavevmode
    \epsfxsize=6in
    \epsfbox{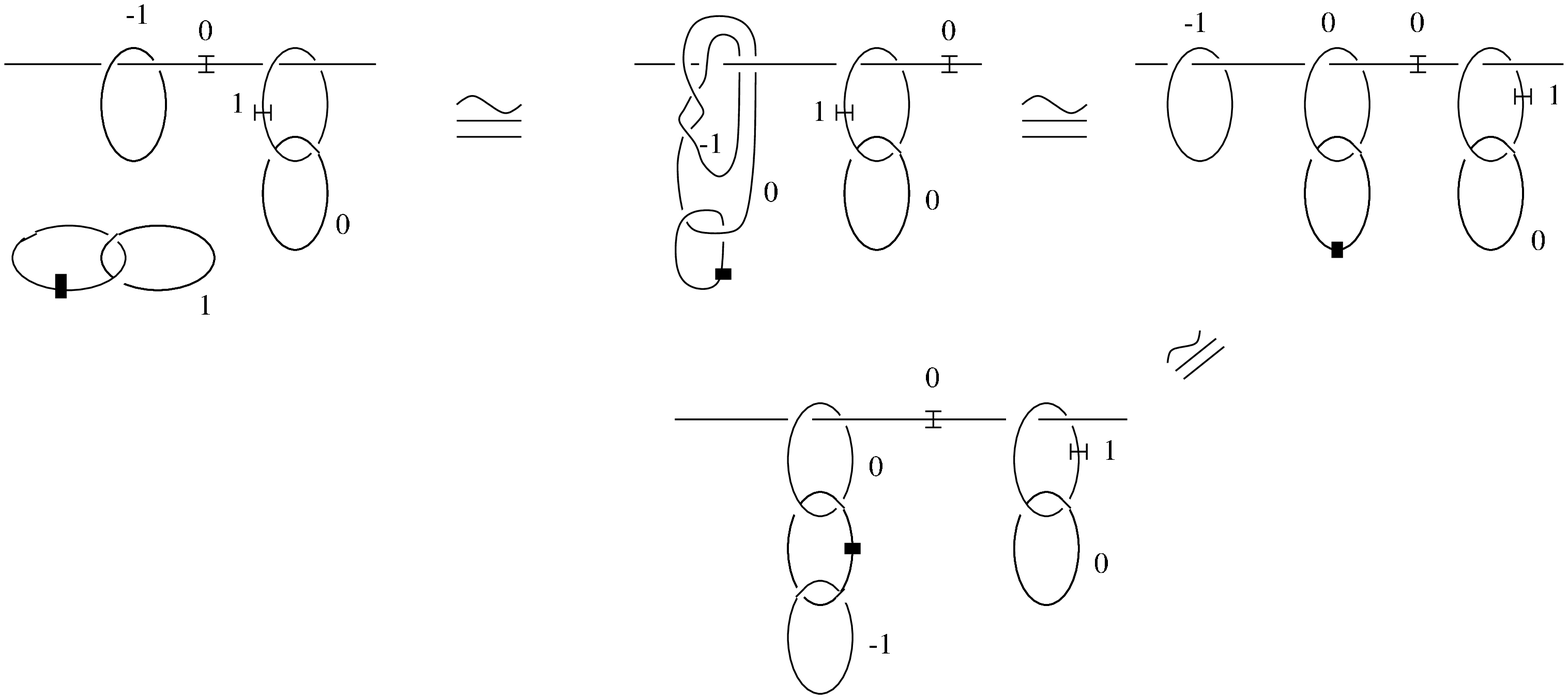}
  \end{center}
   \centerline{\bf Figure 14: Introducing 1-handles.}

\vfill \newpage


\begin{figure}[htbp]
    \leavevmode
    \epsfxsize=6in
    \epsfbox{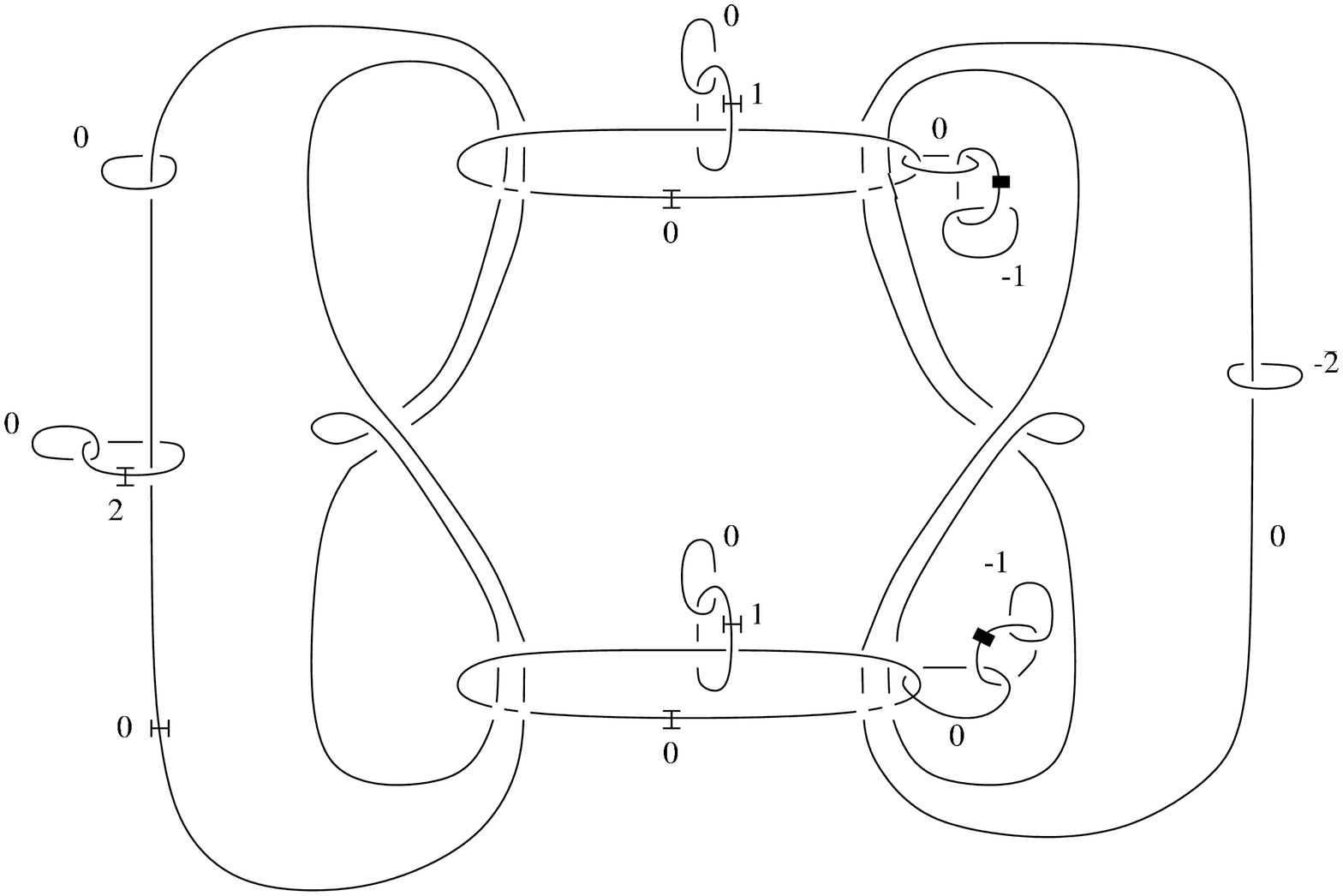} 
   \centerline{\bf Figure 15: $(N_2,T^2)\# (N_2,T^2)\#(S^2\times S^2)$.}
\end{figure}


\nd The handle slides in figure 16 will  make the right side look
exactly like a right-side-up nucleus.
Now, introduce two canceling 2-handle/3-handle pairs. Slide one of
the new 2-handles over the simple $0I$ component, then
use the 2-handles dual to the $1I$ and complicated $0I$
components to arrange the new 2-handle as in figure 17. Repeat with the
second new 2-handle.

\begin{figure}[htbp]
  \begin{center}
    \leavevmode
   \epsfxsize=6in
    \epsfbox{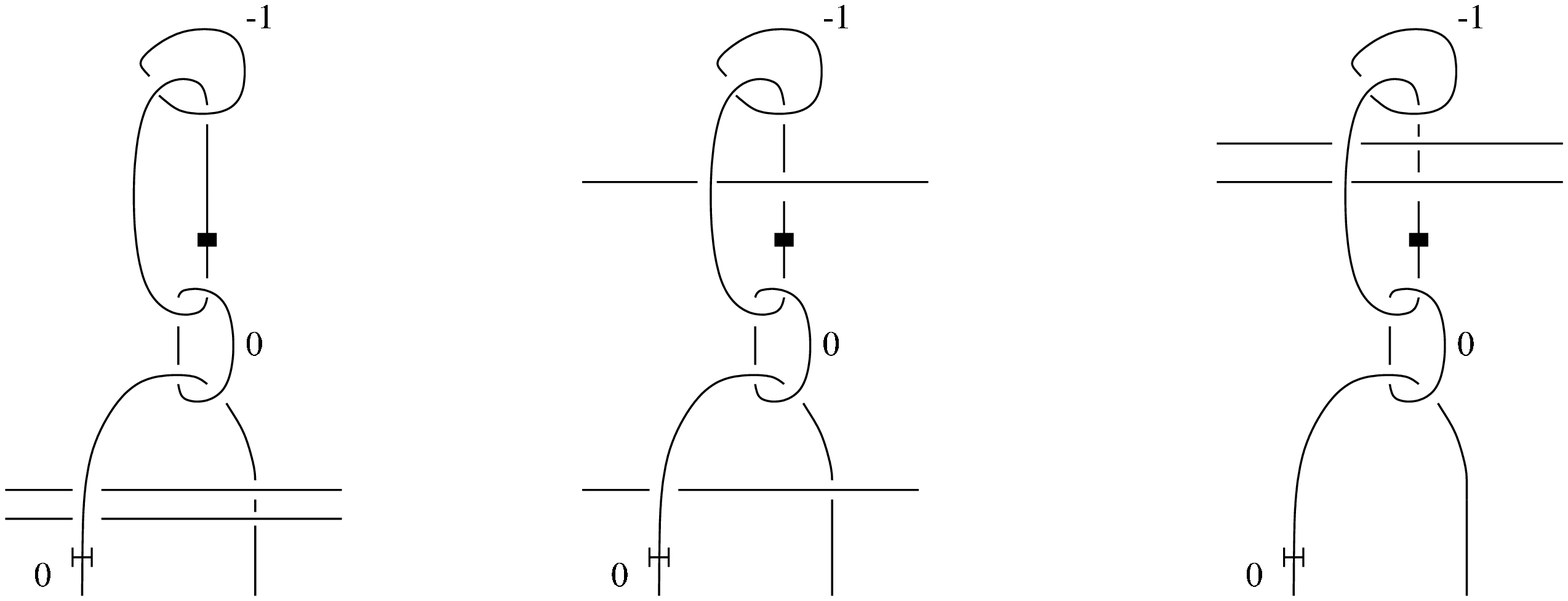}
  \end{center}
   \centerline{\bf Figure 16: Completing a nucleus.}
\end{figure}


  \begin{center}
    \leavevmode
    \epsfxsize=5in
    \epsfbox{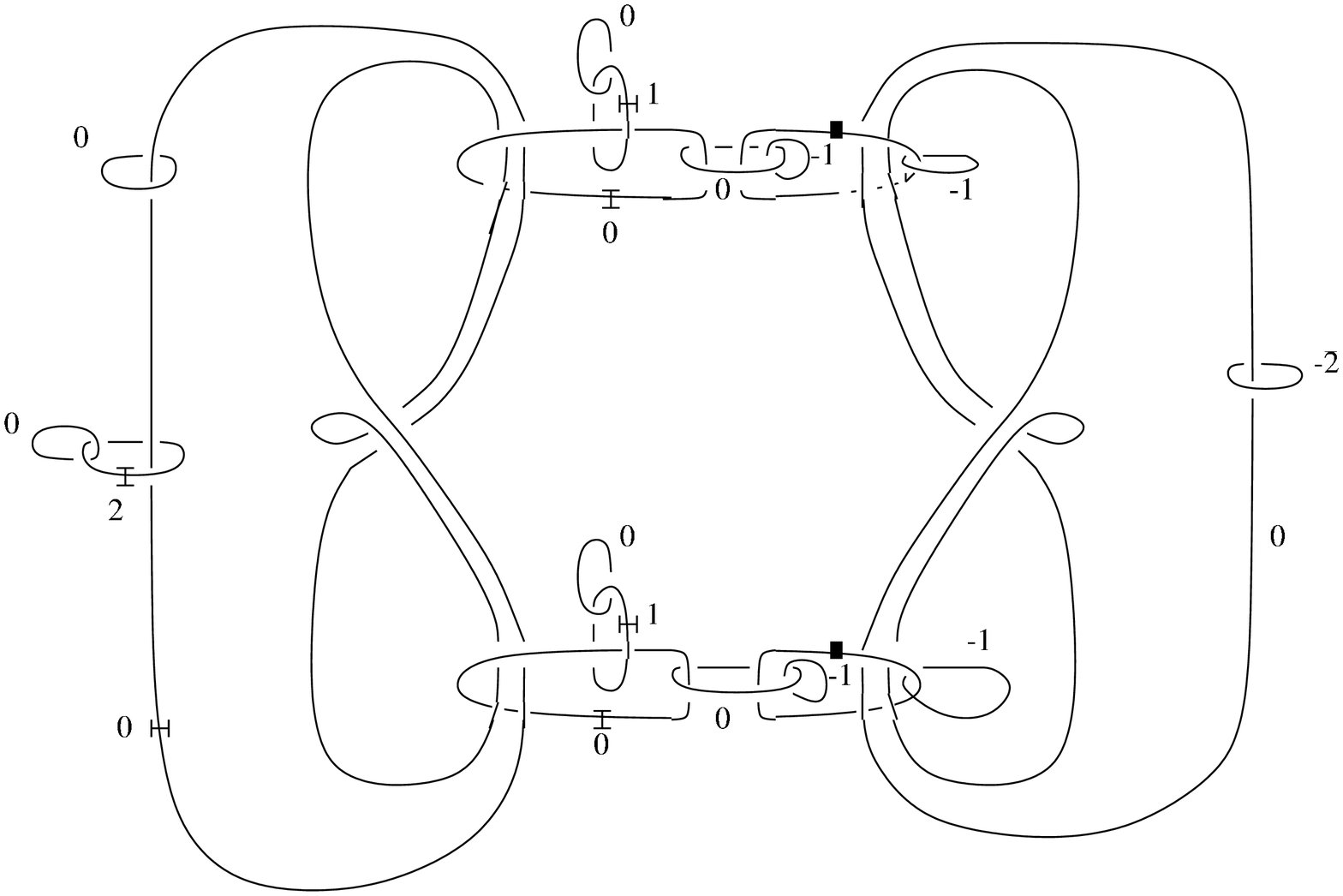}
  \end{center}
   \centerline{\bf Figure 17: $(N_2,T^2)\#(N_2,T^2)\#(S^2\times S^2)$.}

\vskip.1in


\nd Adding the 2-handles in the five dimensional cobordism corresponds
to the handle slides in figure 18. The handle slides in this figure
show that
$(N_2,T^2)\#(N_2,T^2)\#(S^2\times S^2)
  \cong N_2\# N_2\#(S^2\times S^2)\#(S^2\times S^2)$.
This argument may be repeated on each $(N_2,T^2)\#(N_2,T^2)$.
This will show that
$$\begin{array}{rl}
  X_N\#(S^2\~x S^2)
   &\cong (K3)^{\# N-1}\#(S^2\times S^2)^{\# N} \#(\overline{\C P^2})^{\#7}
     \#W_1\cup M\cup E_8  \\
   &\cong (K3)^{\# N-1} \#(S^2\times S^2)^{N-1}\#(S^2\~x S^2)
     \#(\overline{\C P^2})^{\#7}\#W_1\cup M\cup E_8  \\
   &\cong(\overline{\C P^2})^{14N}\#(S^2\times S^2)^{N-1}
     \#(S^2\~x S^2)\#(W_1\cup M\cup W_1)^{\# N}  \\
   &\cong (\C P^2)^{\#4N}\#(\overline{\C P^2})^{\#20N}.
   \end{array} $$

  \begin{center}
    \leavevmode
    \epsfxsize=5.6in
    \epsfbox{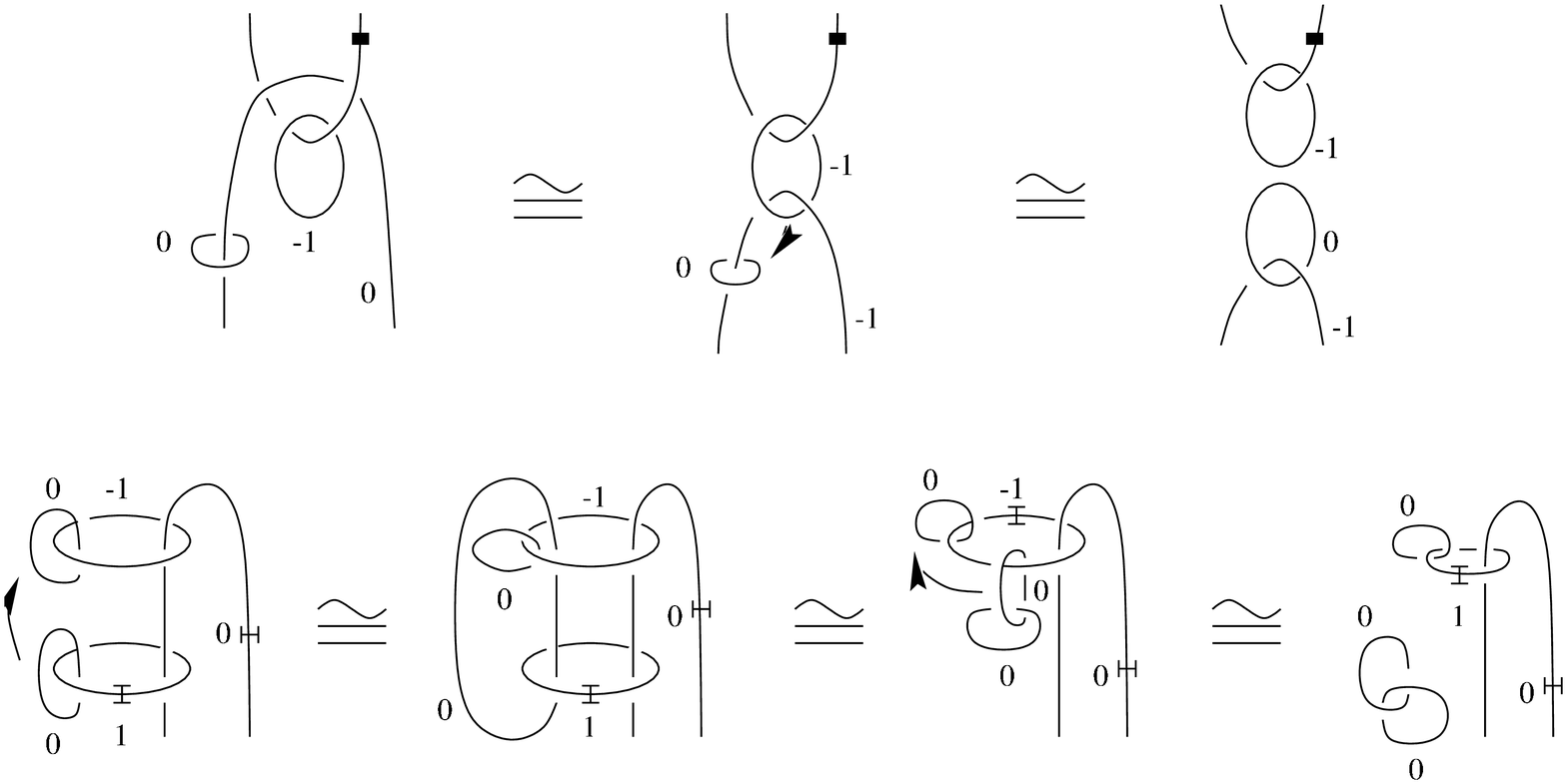}
  \end{center}
   \centerline{\bf Figure 18 $(N_2,T^2)\#(N_2,T^2)\#(S^2\times S^2)
   \cong N^{\#2}_2\#(S^2\times S^2)^{\#2}$.}
  %

\vskip.1in


\nd In the above argument, $M$ is the complement of two $E_8$ manifolds in
K3. Figure 19 displays handle decompositions of $M$ and $W_1\#M\#W_1$.


  \begin{center}
    \leavevmode
    \epsfxsize=5.6in
    \epsfbox{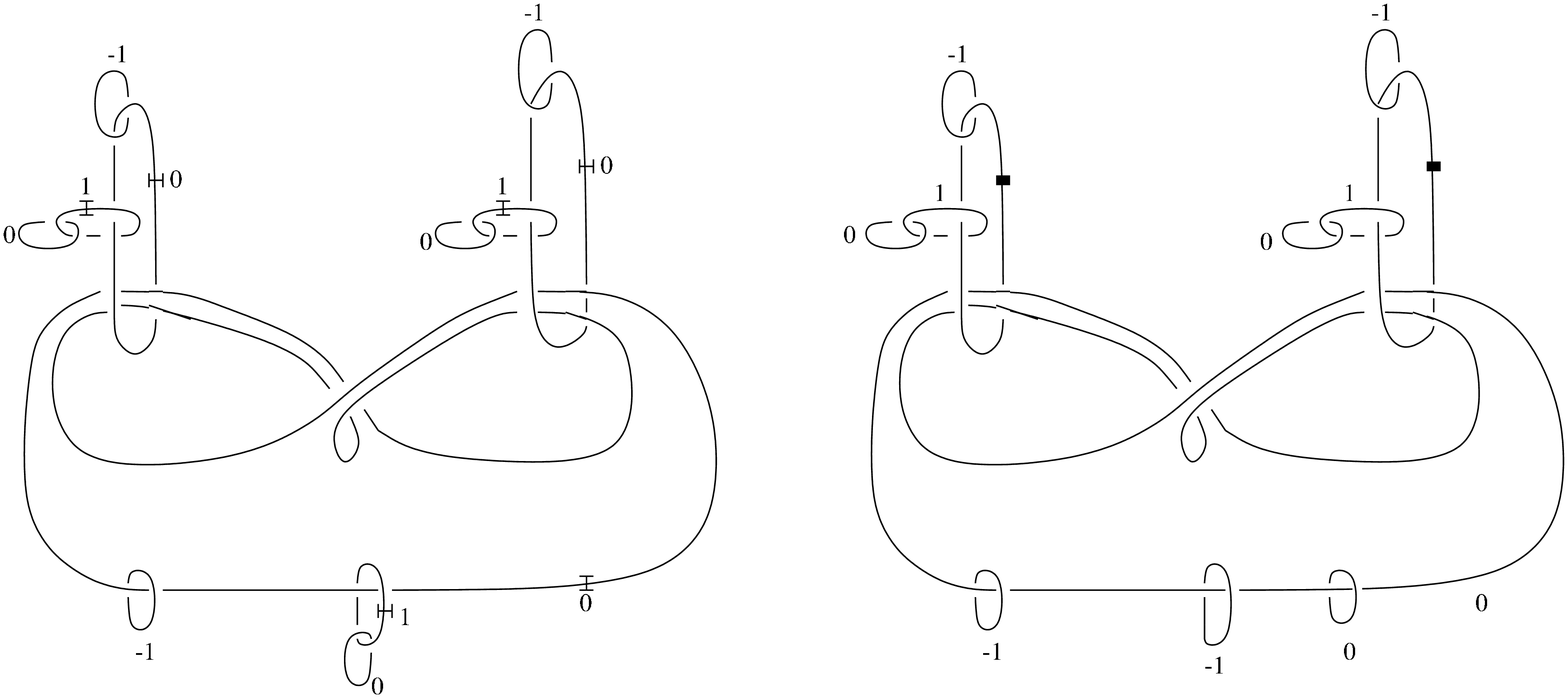}
  \end{center}
   \centerline{\bf Figure 19: $M$ and $W_1\cup M\cup W_1$.}

\vskip.1in


We will now discuss the effect of a single stabilization on a manifold
fiber summed with an $S^3\times S^1$ along a knot cross a circle.
Let $K_1$ and $K_2$ be two knots related by a single crossing change.
By Markov moves, the relevant crossing may be assumed to be in the
lower right corner of a braid representation of $K_1$.
If $S^3-\oN(K_1)$ is described with an extra non-interacting $+1$ Dehn
surgery, then the manifold $(N_2,T)\#(S^3\times S^1, K_1\times S^1)$
will have the handle decomposition displayed in figure 20.
All unlabeled 2-handles are zero framed.

To obtain figure 21, take the connected sum with
$S^2\~x S^2$ and slide handles. Now add two canceling
2-handle/3-handle pairs and one
1-handle/2-handle pair (figure 22). From here a long series of
handle slides will demonstrate that
$(N_2,T)\#(S^3\times S^1, K_1\times S^1)\#
 S^2\~x S^2\cong(N_2,T)\#(S^3\times S^1,K_2\times S^1)\#S^2\~x S^2$
(figures 23-26). The moves from
figure 25 to figure 26 are illustrated in figure 27. The 1-handle
with feet is redrawn, represented by a circle with a dot.
The rightmost strand may be pulled out from the braid by sliding
it over some of the concentric 2-handles.

Finally notice that one can pass from any knot to the unknot by a
series of crossing changes. Call the resulting sequence of knots
$K_1, K_2, \dots, K_n$,  with $K_n$, the unknot. Then
$$ \begin{array}{rl}
  (N_2,t)\#(S^3\times S^1, K_1\times S^1)\#(S^2\~x S^2)
   &\cong (N_2,t)\#(S^3\times S^1, K_2\times S^1)\#(S^2\~x S^2)\dots\\
   &\cong (N_2,t)\#(S^3\times S^1, K_n\times S^1)\#(S^2\~x S^2) \\
   &\cong N_2\#(S^2 \~x S^2).
  \end{array} $$

The previous argument may be distilled to prove that any two manifolds
related by a sequence of special moves become diffeomorphic
after one stabilization. This special move is given in figure 28 which 
displays two different ways to
attach a $T^2\times S^2$ to an $I\times T^3$.

\vfill
\pagebreak~

\vskip .6in
  \begin{center}
    \leavevmode
    \epsfxsize=6in
    \epsfbox{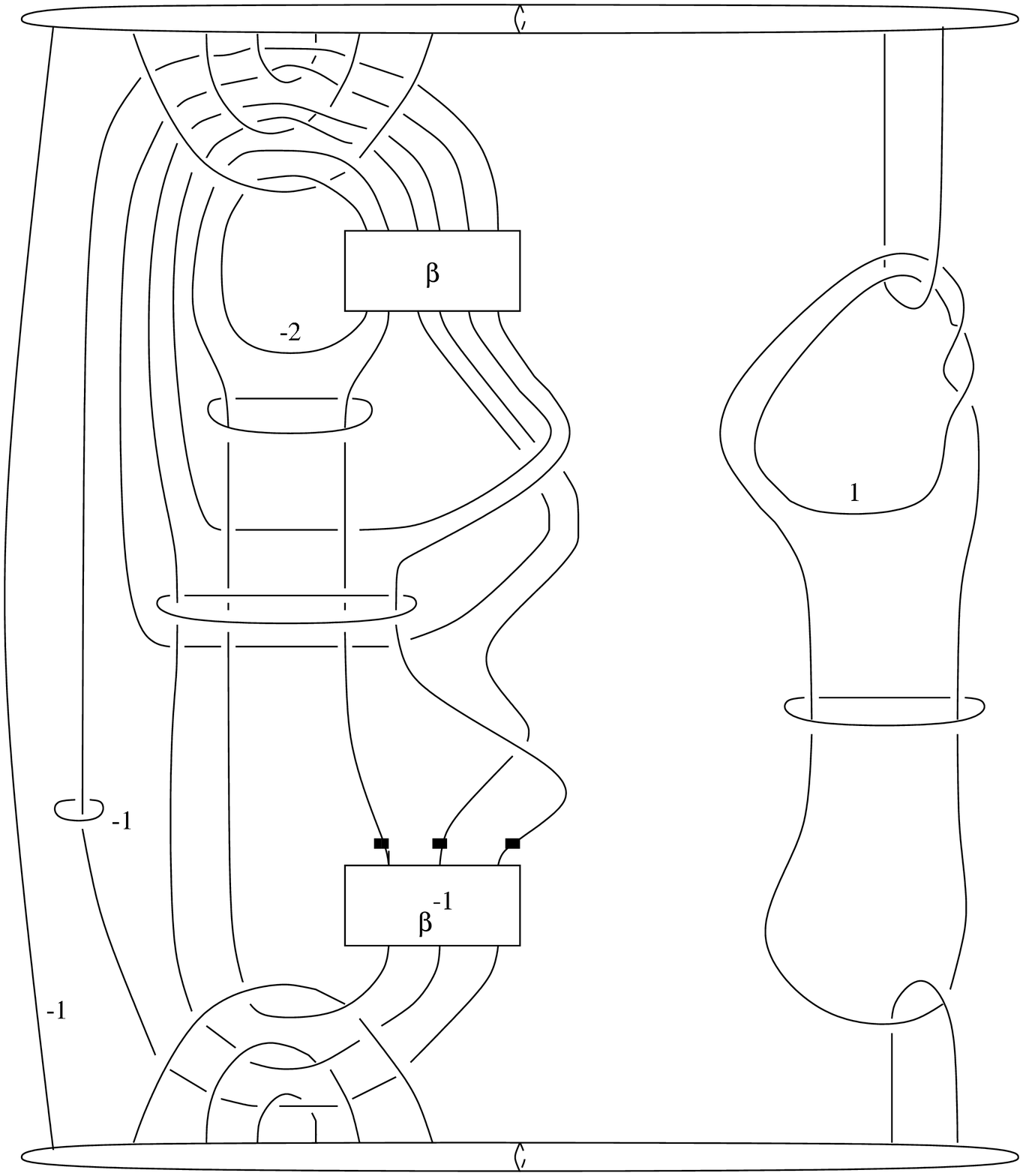}
  \end{center}
   \centerline{\bf Figure 20: $(N_2,T)\#(S^3\times S^1, K_1\times S^1)$.}


\vfill \pagebreak~
\vskip .4in

  \begin{center}
    \leavevmode
    \epsfxsize=6in
    \epsfbox{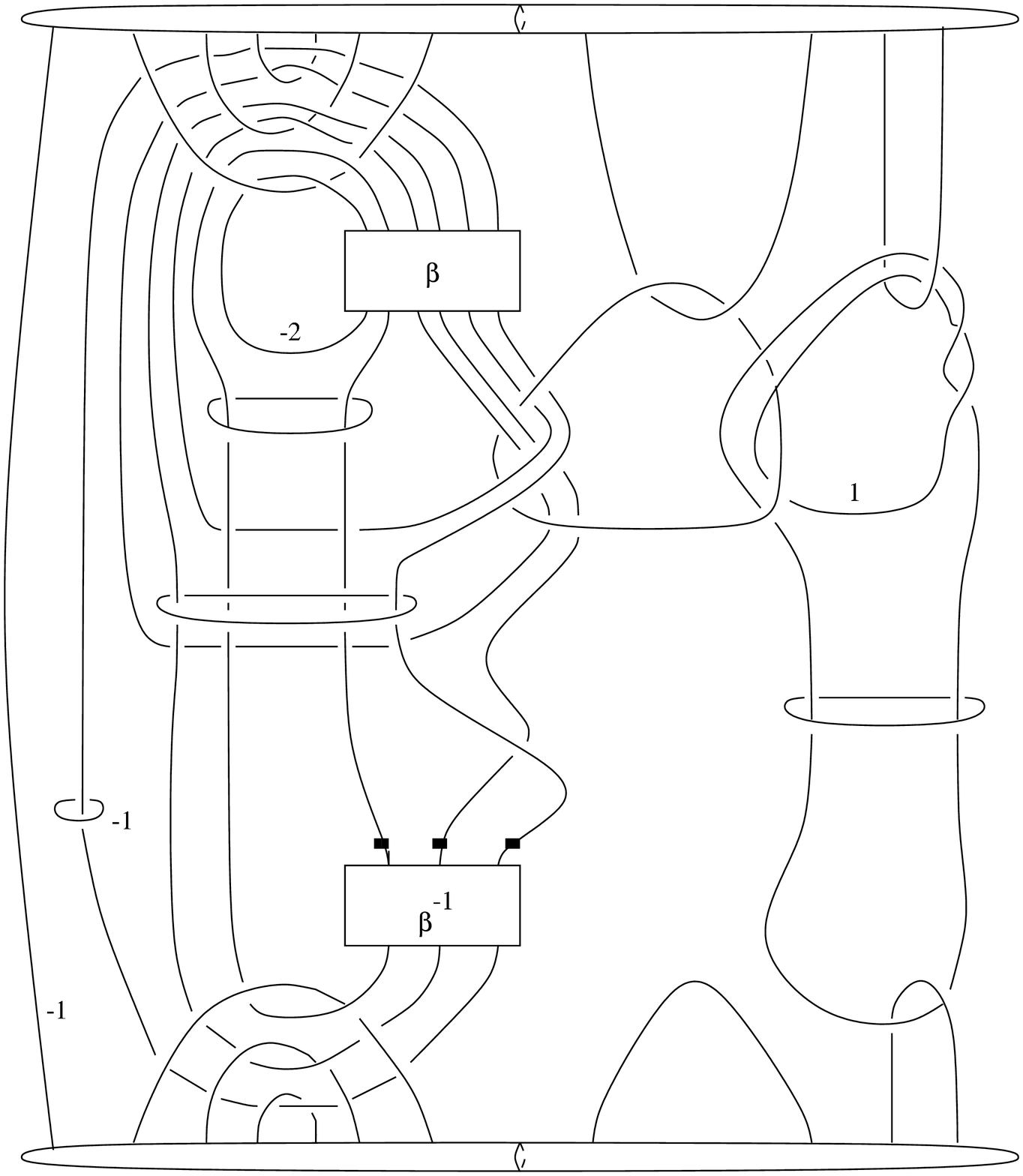}
  \end{center}
   \centerline{\bf Figure 21: 
$(N_2, T)\#(S^3\times S^1, K_1\times S^1)\#(S^2\~x S^2)$.}


\vfill \pagebreak~

  \begin{center}
    \leavevmode
    \epsfxsize=6in
    \epsfbox{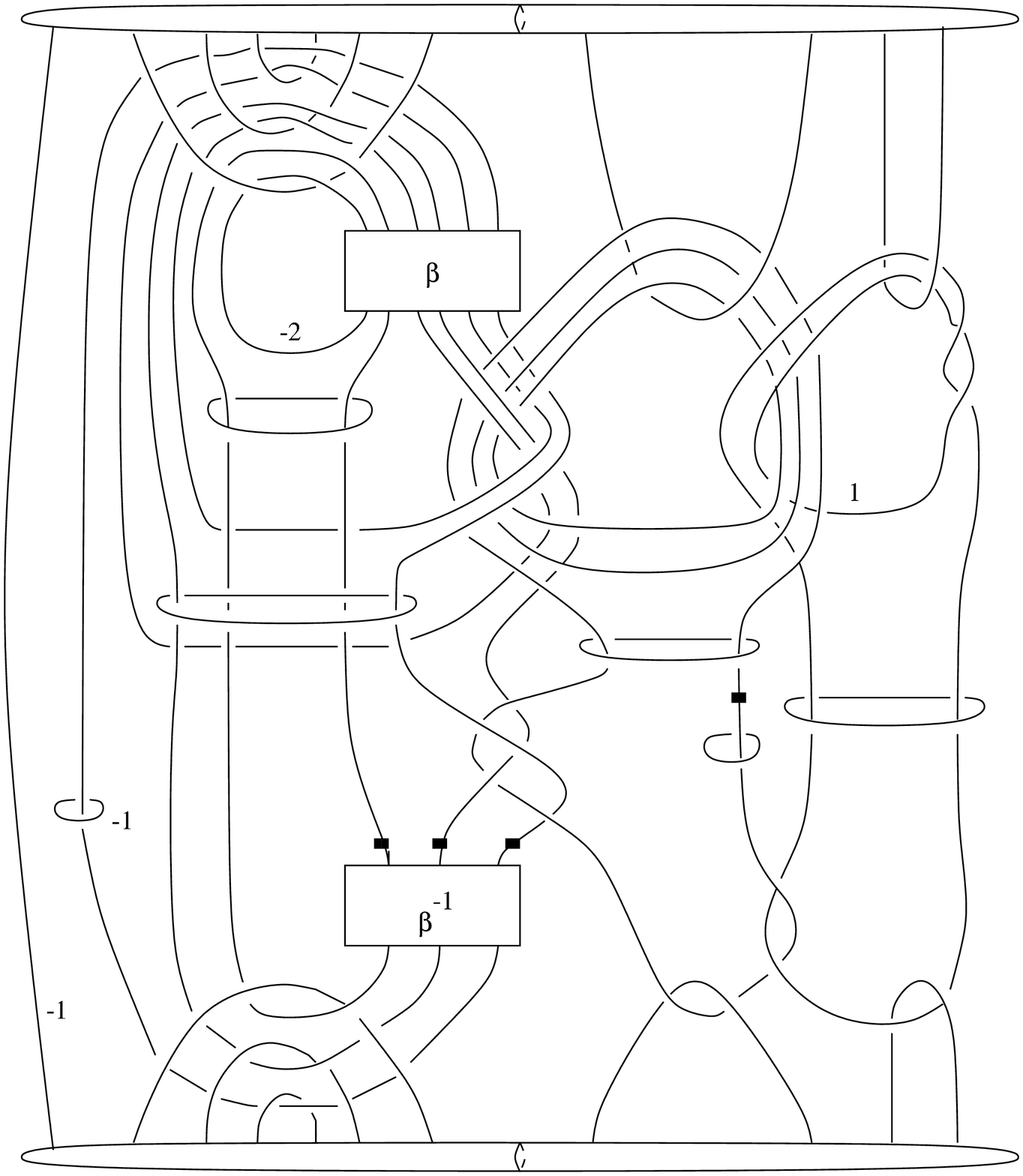}
  \end{center}
   \centerline{\bf Figure 22:
$(N_2,T)\#(S^3\times S^1, K_1\times S^1)\#(S^2 \~x S^2)$.}


\vfill \pagebreak~ \vskip .6in

\begin{figure}[htbp]
  \begin{center}
    \leavevmode
    \epsfxsize=6in
    \epsfbox{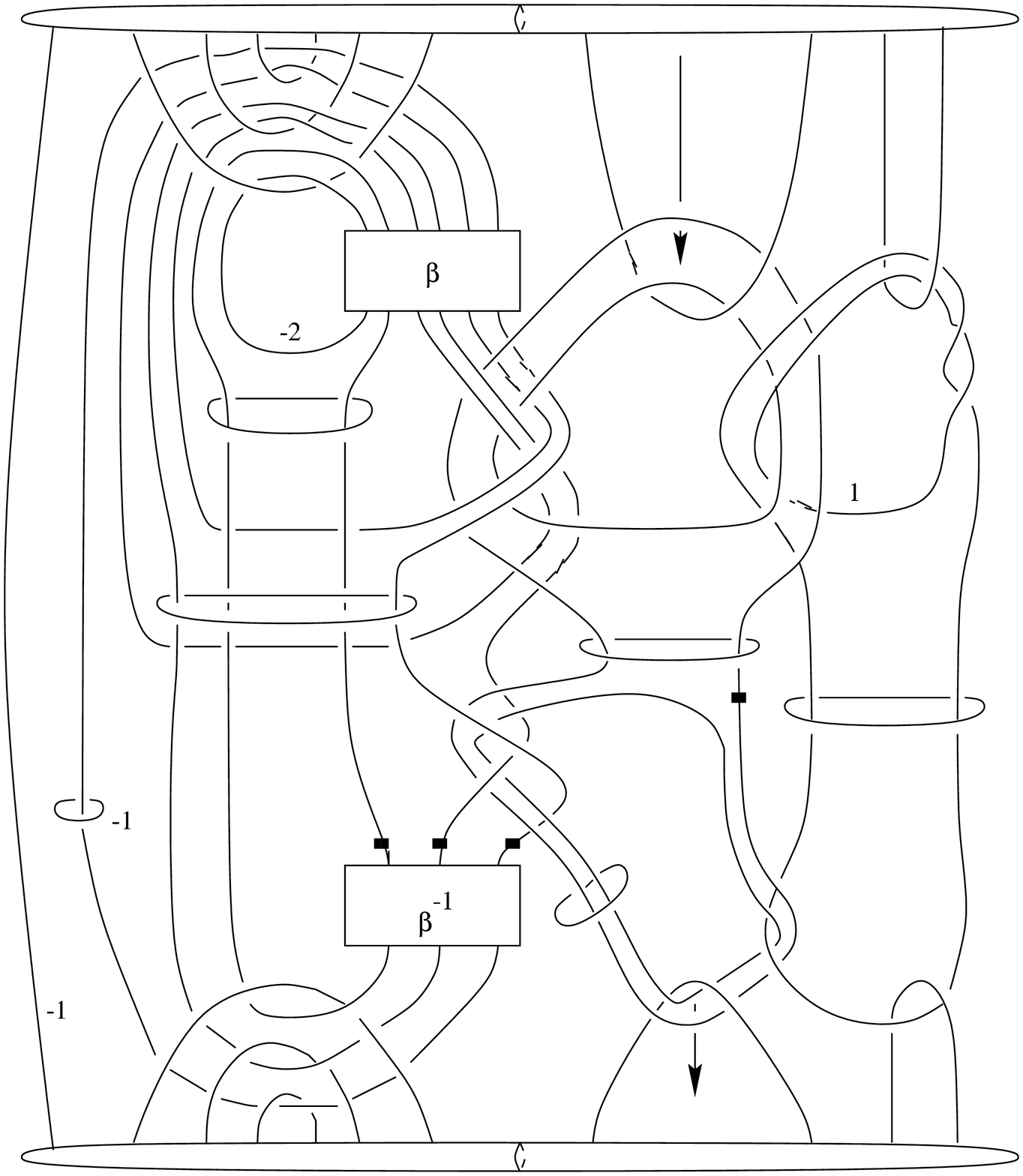}
  \end{center}
   \centerline{\bf Figure 23:
$(N_2,T)\#(S^3\times S^1, K_1\times S^1)\#(S^2\~x S^2)$.}
\end{figure}


\vfill \pagebreak~ \vskip .6in

\begin{figure}[htbp]
  \begin{center}
    \leavevmode
    \epsfxsize=6in
    \epsfbox{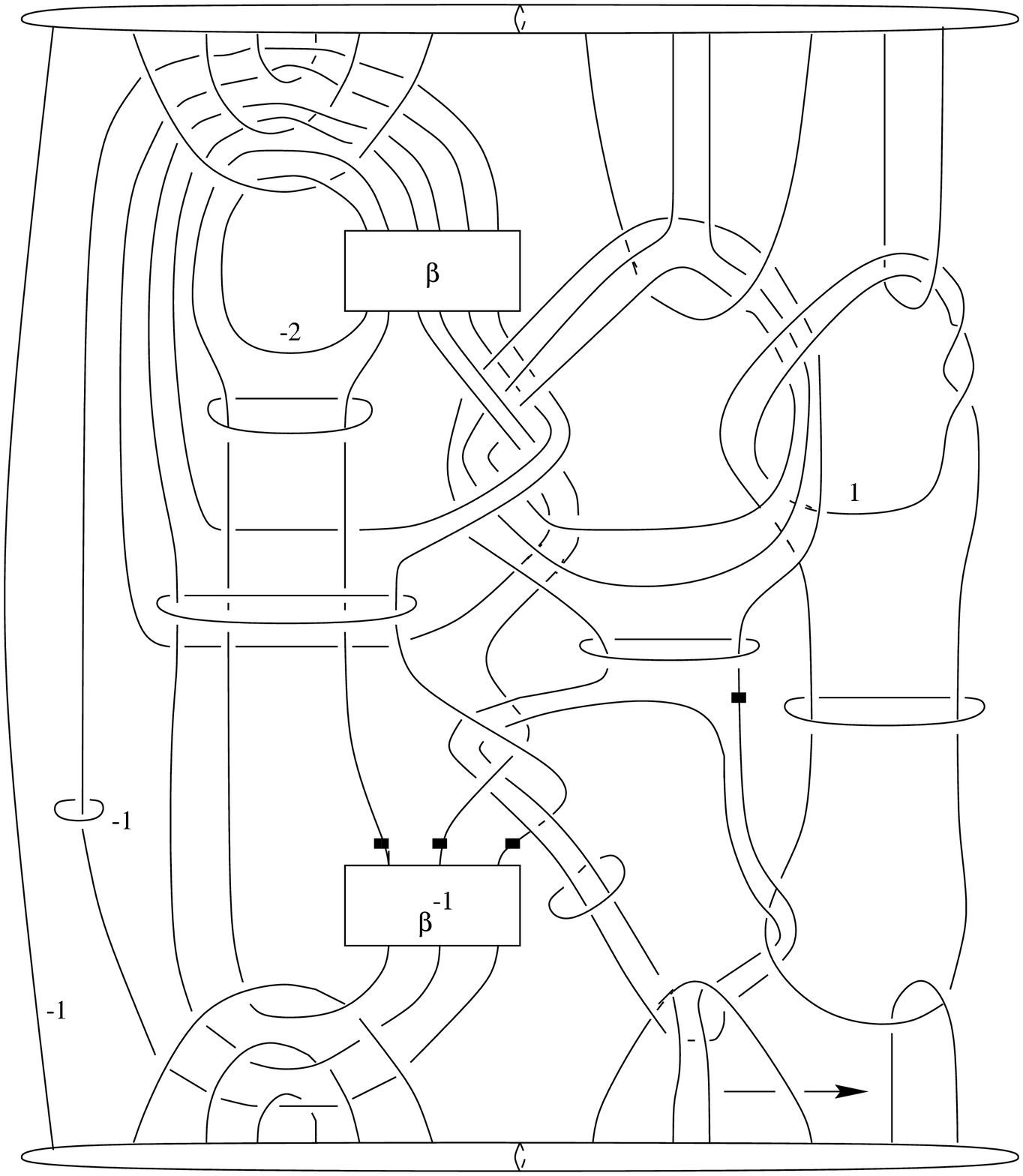}
  \end{center}
   \centerline{\bf Figure 24:
$(N_2,T)\#(S^3\times S^1, K_1\times S^1)\#(S^2\~x S^2)$.}
\end{figure}


\vfill \pagebreak~ \vskip .6in

  \begin{center}
    \leavevmode
    \epsfxsize=6in
    \epsfbox{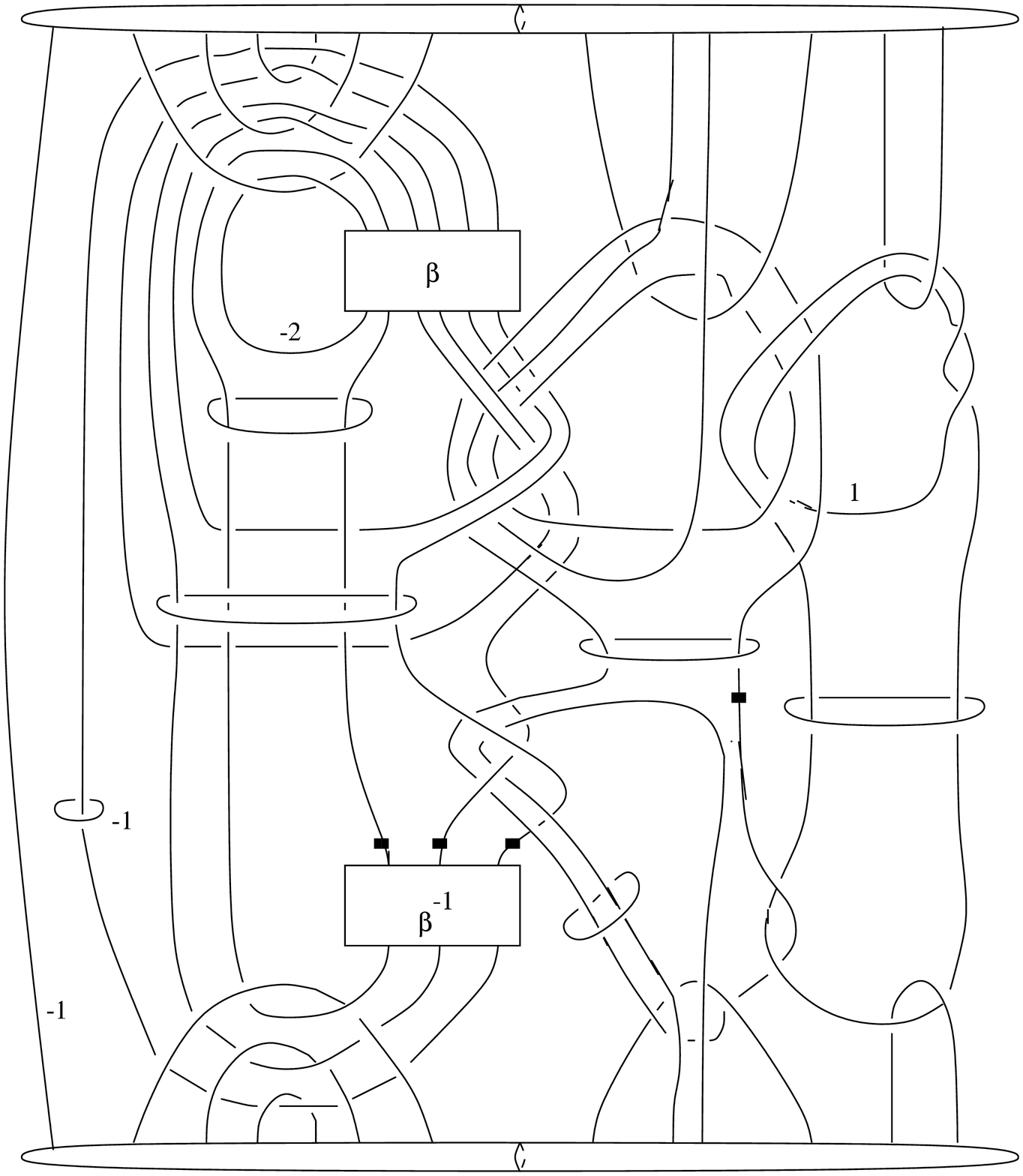}
  \end{center}
   \centerline{\bf Figure 25:
$(N_2,T)\#(S^3\times S^1, K_1\times S^1)\#(S^2\~x S^2)$.}


\vfill \pagebreak~ \vskip .6in
  \begin{center}
    \leavevmode
    \epsfxsize=6in
    \epsfbox{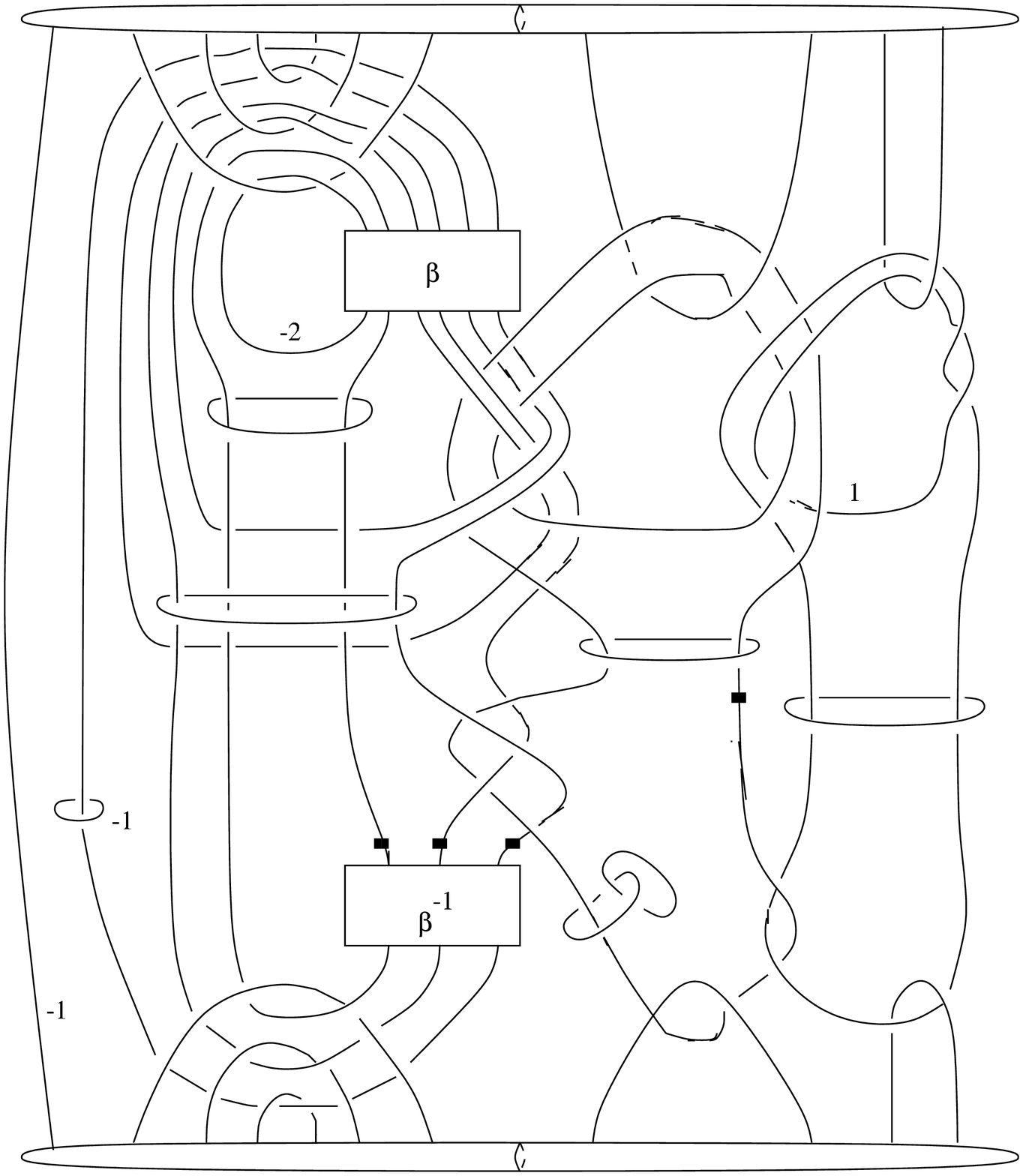}
  \end{center}
   \centerline{\bf Figure 26:
$(N_2,T)\#(S^3\times S^1, K_1\times S^1)\#(S^2\~x S^2)$.}

\vfill \pagebreak

\begin{figure}[htbp]
  \begin{center}
    \leavevmode
    \epsfxsize=3in
    \epsfbox{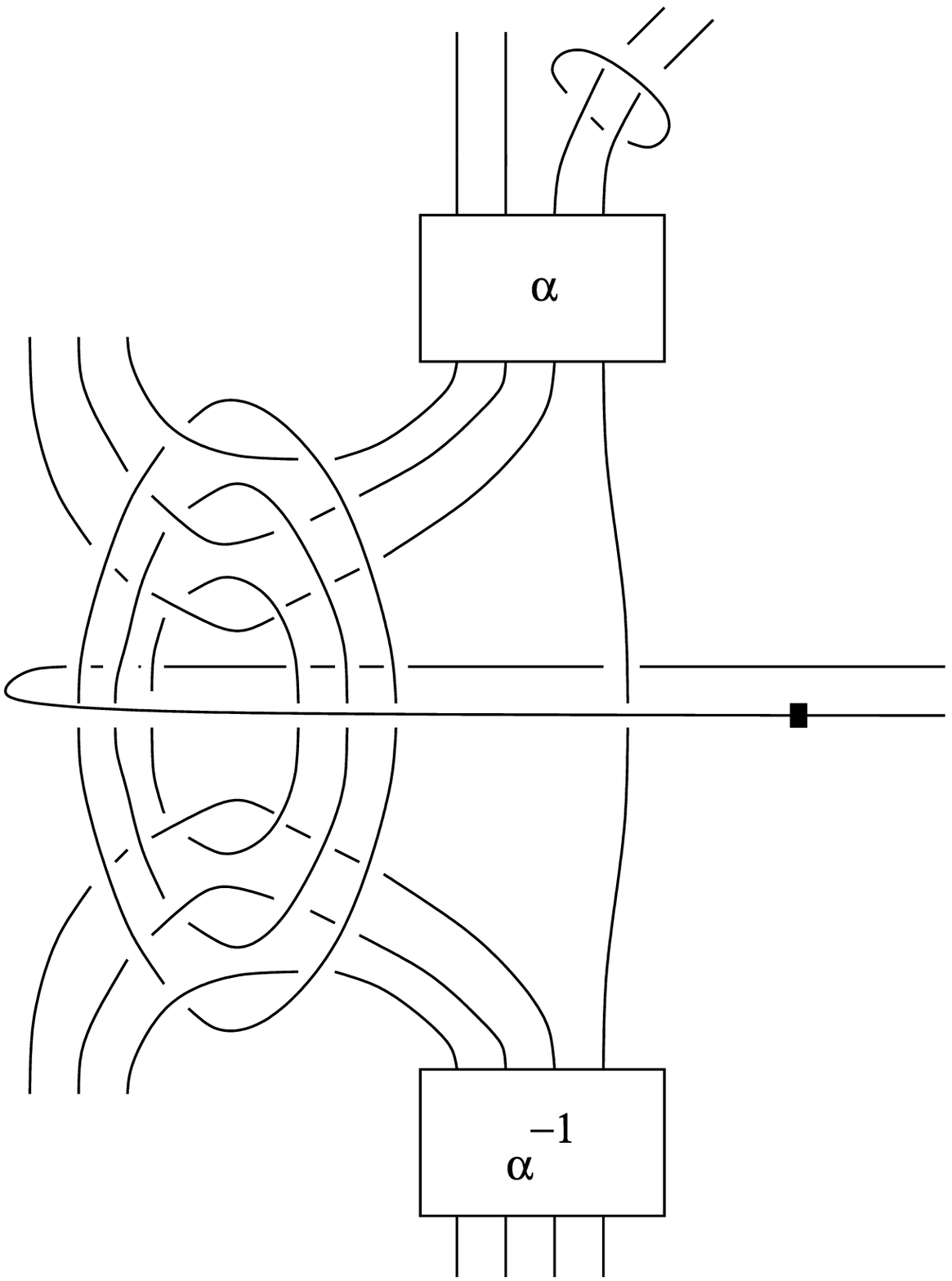}
  \end{center}
   \centerline{\bf Figure 27:
  Pulling a strand away from a braid.}
\end{figure}


  \begin{center}
    \leavevmode
    \epsfxsize=6in
    \epsfbox{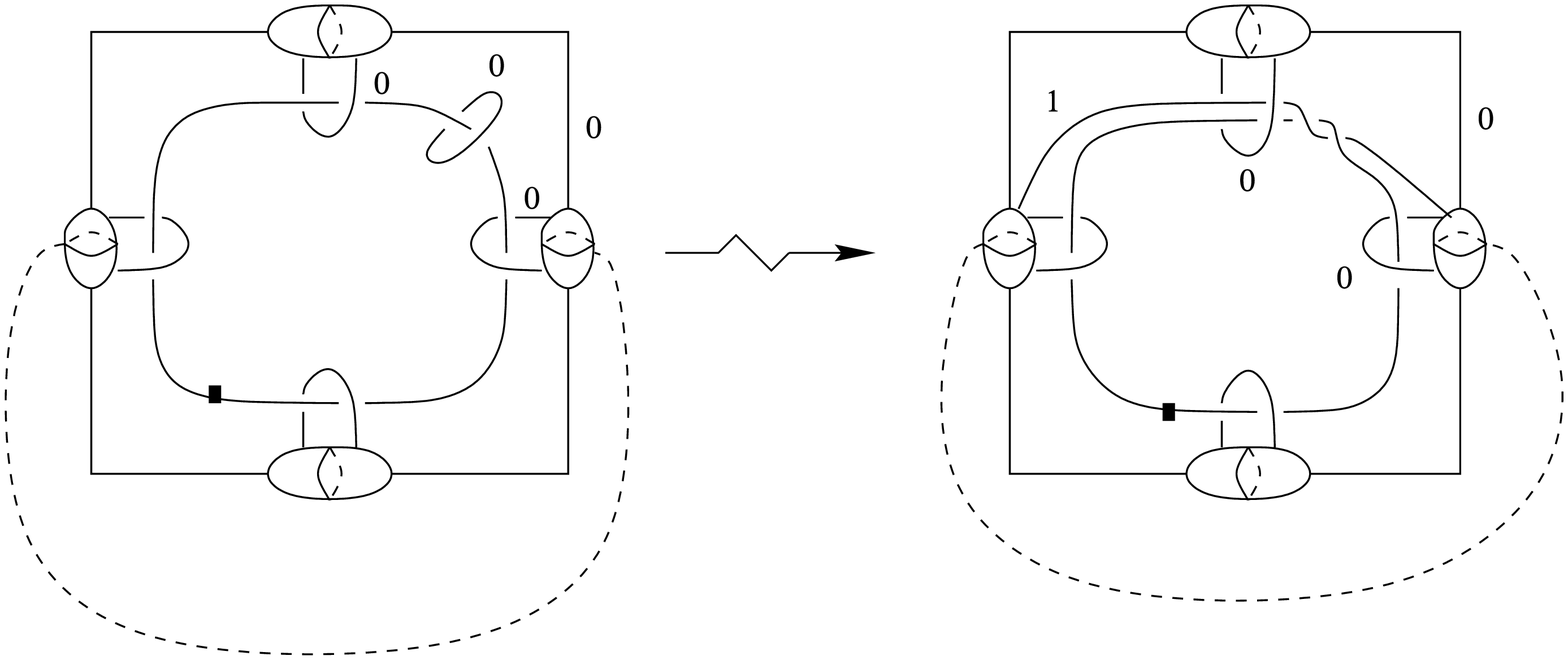}
  \end{center}
   \centerline{\bf Figure 28:
  Geometrically null $+1$ log transform.}


\nd If the dotted line bounds an evenly framed disk in some four-manifold,
we will call the process of cutting out a $T^2\times D^2$ and regluing it
a geometrically null $+1$ log transform. This is just
the product of $+1$ surgery with a circle.
The Kirby calculus in figures
29 and 30 demonstrates the following theorem.

\begin{theorem}
Two manifolds related by a geometrically null
$+1$ log transform become diffeomorphic after one stabilization.
\end{theorem}

\begin{figure}[htbp]
  \begin{center}
    \leavevmode
    \epsfxsize=6in
    \epsfbox{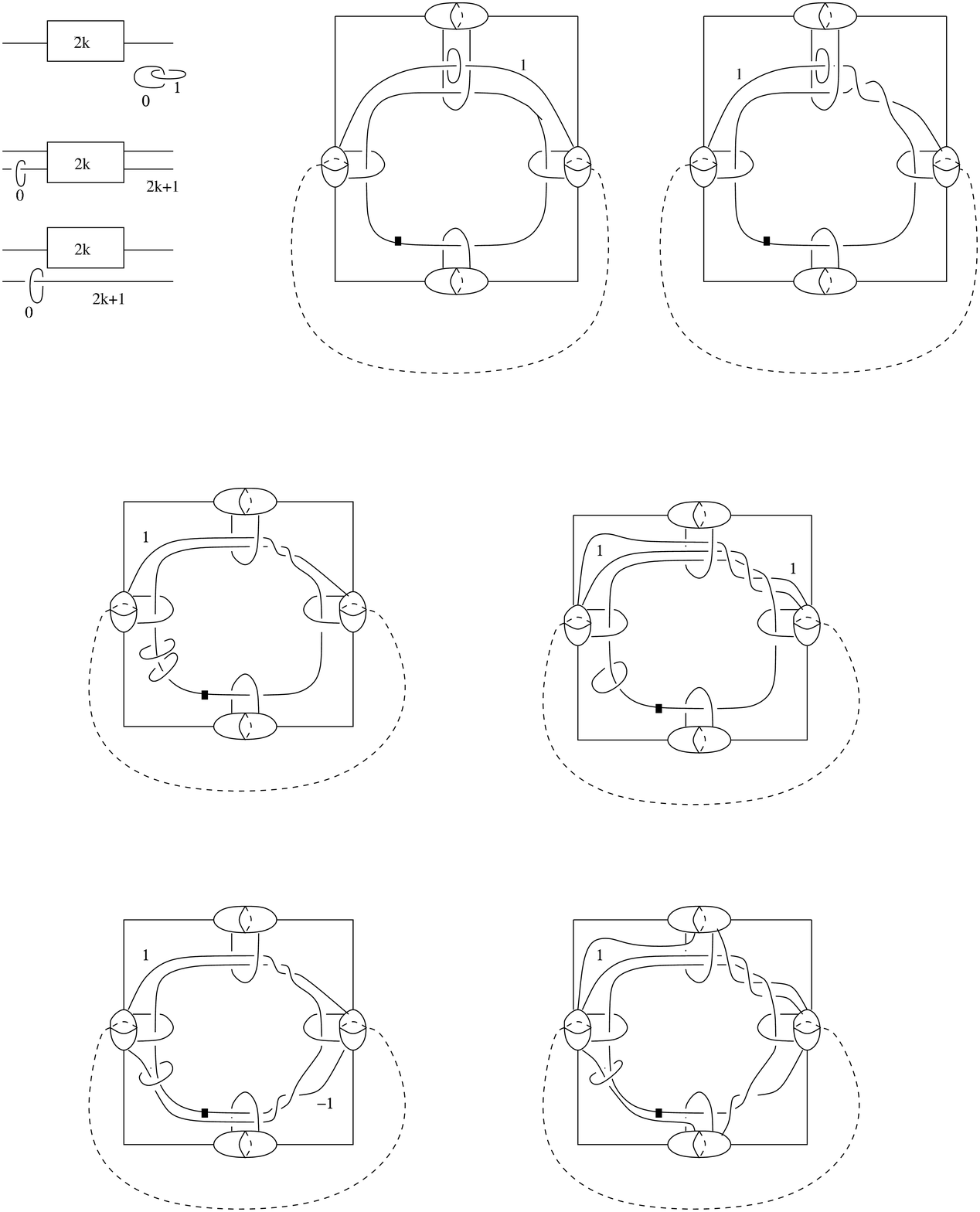}
  \end{center}
   \centerline{\bf Figure 29: Stabilizing a log transform}
\end{figure}


\begin{figure}[htbp]
  \begin{center}
    \leavevmode
    \epsfxsize=6in
    \epsfbox{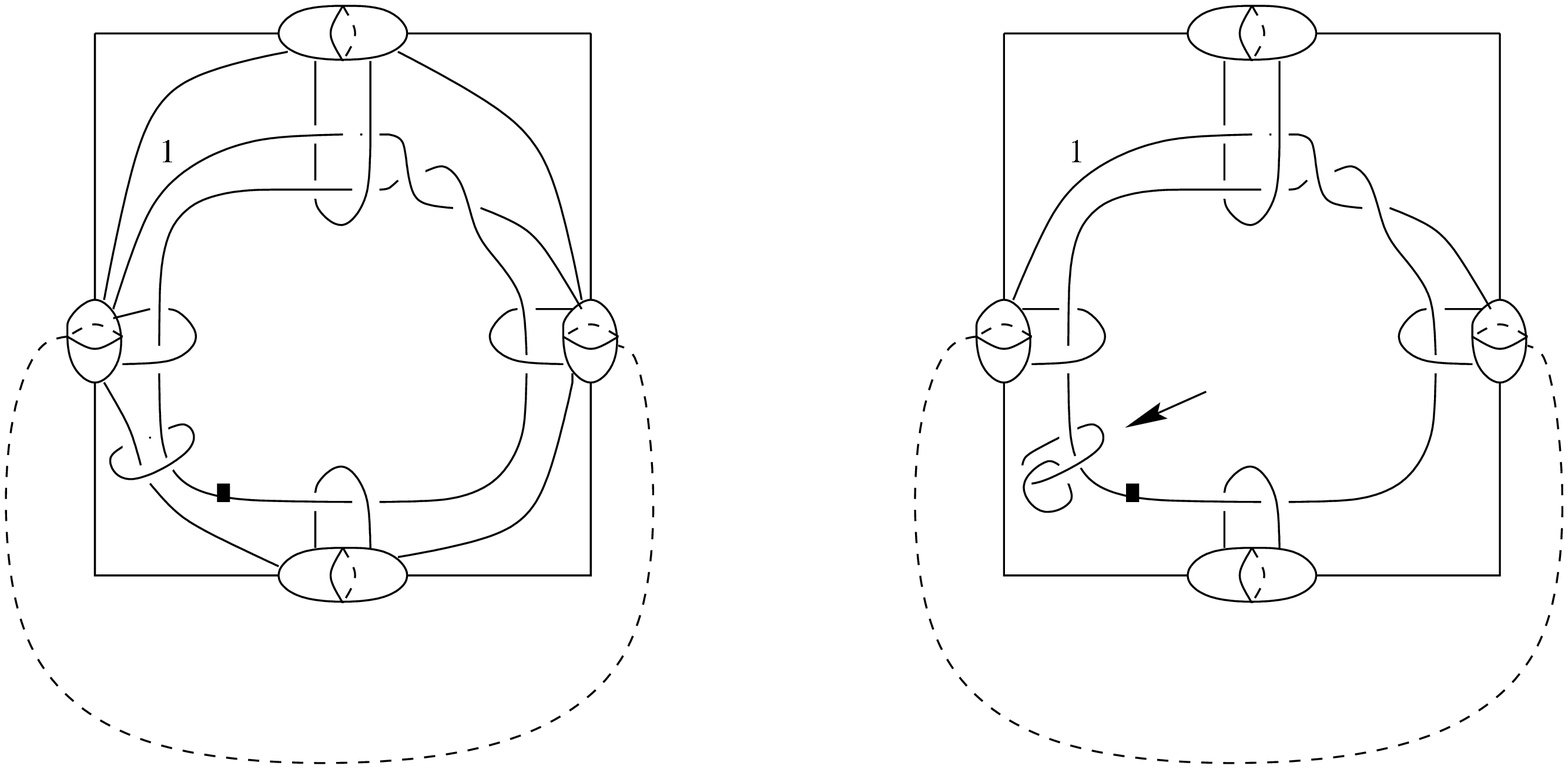}
  \end{center}
   \centerline{
Finish by sliding the labeled handle over the $+1$ framed handle
and reversing}
\centerline{ the moves from the beginning.}
\centerline{\bf Figure 30: Stabilizing a log transform}
\end{figure}



\section*{Seiberg-Witten invariants}

Recall that the Seiberg-Witten series of a smooth 4-manifold with
homology orientation is
$$ SW_X=a_0+\sum aj
  \left( \exp(K_j)+(-1)^{\frac{\chi(x)+\alpha(x)}{4}}
   \exp (-k_j) \right) $$
where the set of basic classes is
$\{\pm K_1, \pm K_2,\dots \pm K_n\}\subseteq H^2(X;\Z)$,
$a_0=SW_X(0)$, and $a_j=SW_X(K_j).$
If $b^+_2(X)>0$,
then $SW_{X\#(S^2\~x S^2)}=0$.
Thus the Seiberg-Witten invariant cannot distinguish the two
manifolds, $X\#(S^2\~x S^2)$ and $Y\#(S^2\~x S^2)$.
One might hope that a diffeomorphism between
 $X\#(S^2\~x S^2)$ and $Y\#(S^2\~x S^2)$
would imply some restriction on the relationship between the
Seiberg-Witten series, $SW_X$ and $SW_Y$. We will compute
the Seiberg-Witten series of all of the manifolds
considered in the previous section. The number of
basic classes, the  rank of the space spanned by the basic classes,
and the coefficients of the Seiberg-Witten series will vary
arbitrarily in each family, $F_N$, of manifolds.

To compute the Seiberg-Witten series, we will use several gluing
formula worked out by Morgan, Mrowka, and Szabo, and utilized
by Fintushel and Stern [17], [16], [5].

\begin{description}
\item[\bf Fact 1:] $SW_{K3}=1$
\item[\bf Fact 2:] $SW_{(X,T)\#(Y,S)}
    =SW_X\cdot SW_Y\cdot(\exp(T)-\exp(-T))^2.$

\item[\bf Fact 3:]
If $\pi_1(X)=1$, $\pi_1(X-T)=1$, $[T]\not=0$ in $H_2(X)$
and $[T]^2=0$,
\item[\qquad\quad] then $SW_{(X,T)\#(S^3\times S^1, K\times S^1)}
   =SW_X\cdot\Delta_K(\exp(2T))$.
\item[\qquad\quad] Here, $\Delta_K$ is the Alexander polynomial of $K$.

\end{description}

The first fact is due to Witten, and is by now well known [21], [19].
The second fact has not yet appeared in the literature, but it is similar
to  the results in [17] and [16].
We have not included the technical hypothesis for the second fact.
The third fact is proved in [5].

Refine our original notation, to denote the tori in $X_N$ by $T_{\alpha,i}$,
with $\alpha=1,2\dots,N$ and $i=1,2,3$ so that
$T_{\alpha,3}=T_{\alpha+1,1}$ for $\alpha=1,\dots,N-1$.
Using this notation, the Seiberg-Witten series of $X_N$ is
$$ SW_{X_N}=\prod^{N-1}_{\alpha=1}
  (\exp(T_{\alpha,3})-\exp(-T_{\alpha,3})). $$
Finally, let
$$\begin{array}{rl}
   & Y_0= X_N, Y_{\alpha+1}=(Y_\alpha,T_{\alpha,2})
         \#(S^3\times S^1, K_{\alpha,2}\times S^1) \\
   & Y^\prime =(Y_N,T_{1,1})\#(S^3\times S^1, K_{1,1}\times S^1)\\
   \hbox{and\quad }&Y=(Y^\prime,T_{N,3})\#(S^3\times S^1, K_{N,3}).
   \end{array}$$
Then the Seiberg-Witten series of $Y$ is
$$ \begin{array}{rl}
  SW_Y=\ds\prod^{N-1}_{\alpha=1}
  & \left[ (\exp(T_{\alpha,3})-\exp(-T_{\alpha,3}))
    \cdot\Delta_{K_{\alpha,2}} (\exp(2T_{\alpha,2}))\right]\\
  & \cdot \Delta_{K_{1,1}} (\exp(2T_{1,1}))
    \cdot \Delta_{K_{N,2}} (\exp(2T{N,2}))
    \cdot\Delta_{K_{N,3}} (\exp(2T_{N,3})).
   \end{array} $$


\end{document}